\documentclass[12pt]{amsart}
\usepackage{amssymb,needspace}
\usepackage[hidelinks]{hyperref}

\topmargin \dimexpr-0.5cm-8pt\relax
\newtheorem{thm}{Theorem}[subsection]
\makeatletter
\@addtoreset{thm}{section}
\makeatother

\newtheorem{lem}[thm]{Lemma}
\newtheorem{cor}[thm]{Corollary}
\newtheorem{prop}[thm]{Proposition}
\newtheorem{ex}[thm]{Example}
\newtheorem{claim}[thm]{Claim}
\theoremstyle{definition}
\newtheorem{defi}[thm]{Definition}
\theoremstyle{remark}
\newtheorem{rem}[thm]{Remark}
\numberwithin{equation}{section}

\newcommand{\C}{\mathbb{C}}
\newcommand{\F}{\mathbb{F}}
\newcommand{\Q}{\mathbb{Q}}
\newcommand{\Z}{\mathbb{Z}}
\DeclareMathOperator{\Cong}{Cong}
\DeclareMathOperator{\Discr}{Discr}
\DeclareMathOperator{\Prod}{Prod}
\DeclareMathOperator{\RRes}{RRes}
\DeclareMathOperator{\Res}{Res}
\DeclareMathOperator{\Kh}{Kh}
\DeclareMathOperator{\Gal}{Gal}
\DeclareMathOperator{\End}{End}
\DeclareMathOperator{\GL}{GL}
\DeclareMathOperator{\Sym}{Sym}
\DeclareMathOperator{\dl}{dl}

\begin{document}

\title[Automorphisms satisfying an irreducible identity]
  {Finite groups with a fixed-point-free automorphism
   satisfying an irreducible identity}

\author[D. Burde]{Dietrich Burde}
\address{Fakult\"at f\"ur Mathematik\\
Universit\"at Wien\\
Oskar-Morgenstern-Platz 1\\
1090 Wien\\
Austria}
\email{dietrich.burde@univie.ac.at}

\author[W. A. Moens]{Wolfgang A. Moens}
\address{Fakult\"at f\"ur Mathematik\\
Universit\"at Wien\\
Oskar-Morgenstern-Platz 1\\
1090 Wien\\
Austria}
\email{wolfgang.moens@univie.ac.at}

\date{}

\begin{abstract}
For any given irreducible polynomial $f(x)=a_0+a_1\cdot x+\cdots+a_d\cdot x^d\in\Z[x]$,
we study the structure of finite groups $(G,\cdot)$ that have a fixed-point-free-automorphism
$\alpha$ satisfying the condition
\[
  \{g^{a_0}\cdot\alpha(g)^{a_1}\cdots\alpha^d(g)^{a_d}\mid g\in G\}=\{1\}.
\]
We show that such a group $G$ is nilpotent of class at most $d^{2^d}$, provided that
$|G|$ is coprime to an explicit, non-zero, integer invariant $\iota(f(x))$ of $f(x)$.
We do this by extending the solutions of the Frobenius conjecture. We also extend the
solution of the restricted Burnside problem in prime exponent and the classification of
the $n$-abelian groups. The analogous set-up for reducible polynomials $f(x)$ will be
treated in two subsequent papers.
\end{abstract}

\maketitle

\section{Introduction}\label{sec:1}

We study the structure of finite groups $(G,\cdot)$ having an automorphism $\alpha$ that is fixed-point-free,
by which we mean that $\alpha$ fixes only the neutral element $1$ of $G$. Such a group is solvable according
to the classification of the finite simple groups \cite{Rowley}. Whether the derived length of $G$ can be
bounded by a function of $|\alpha|$ is still an open question. But if the automorphism has \emph{prime} order
$p$, then group is even known to be nilpotent by a theorem of Thompson \cite{Thompson1959a}. Higman
\cite{Higman} showed that the class of $G$ can then be bounded from above by some function of $p$ and
Kreknin--Kostrikin \cite{KrekninKostrikin} later found the explicit upper bound $(p-1)^{2^{(p-1)}}$. The
latter three results settled the so-called Frobenius conjecture. We now recall that a fixed-point-free
automorphism of prime order $p$ satisfies the ``polynomial identity''
\[
\{g\cdot \alpha(g)\cdot \alpha^2(g)\cdots \alpha^{p-1}(g)\mid g\in G\}=\{1\}
\]
and the coprimeness condition $\gcd(|G|,p)=1$. In this paper, we will generalize the above results from the
cyclotomic polynomials $\Phi_p(x)=1+x+x^2+\cdots+x^{p-1}$ to all irreducible polynomials in the following
fashion.

\begin{cor}\label{1.0.1}
For every irreducible polynomial $f(x)=a_0+a_1\cdot x+\cdots+a_d\cdot x^d\in \Z[x]$, there is a non-zero
integer $\iota(f(x))$ with the following property. Let $G$ be a finite group with a fixed-point-free
automorphism $\alpha$ such that
\[
\{g^{a_0}\cdot \alpha(g)^{a_1}\cdot \alpha^2(g)^{a_2}\cdots \alpha^d(g)^{a_d}\mid g\in G\}=\{1\}
\]
and $\gcd(|G|,\iota(f(x)))=1$. Then $G$ is nilpotent of class at most $d^{2^d}$.
\end{cor}

The constant $\iota(f(x))$ will be made explicit in Definition~\ref{8.1.1}. We will also obtain substantially
sharper bounds on the class of $G$ in Corollary~\ref{8.1.4} and we will obtain an optimal bound for
``almost-all'' polynomials $f(x)$ in Remark~\ref{8.1.5}. Examples of Gross \cite{Gross} show moreover that
neither the irreducibility of $f(x)$ nor the coprimeness condition $\gcd(|G|,\iota(f(x)))=1$ can be removed.

Finite groups having a fixed-point-free automorphism of \emph{composite} order have also been studied
extensively in the literature. In two subsequent papers, we will correspondingly study the structure of finite
groups having a fixed-point-free automorphusms satisfying a \emph{reducible} ``polynomial identity'' $f(x)$.
Cf. Remark~\ref{8.1.7}.

Our two main theorems are more general and detailed than the corollary. In order to state them efficiently, we
introduce identities of endomorphisms.

\subsection{Identities of endomorphisms}\label{sec:1.1}

\begin{defi}\label{1.1.1}
Let $G$ be a group with endomorphism $\gamma$. The polynomial $f(x)= a_0+a_1\cdot x+\cdots+a_d\cdot x^d\in
\Z[x]$ is an \emph{ordered identity} of the endomorphism if the map $f(\gamma):G\longrightarrow G:g\longmapsto
g^{a_0}\cdot \gamma(g)^{a_1}\cdots \gamma^d(g)^{a_d}$ vanishes.
\end{defi}

These ordered identities naturally generalize various notions in the literature.

\begin{ex}\label{1.1.2}
Let $n\in \Z_{\geqslant 1}$.
\begin{itemize}
\item[(a)]
If the exponent of a group $G$ divides $n$, then the constant polynomial $f(x):=n$ is an 
ordered identity of the endomorphism $x\mapsto x$. And vice-versa. Finite groups of 
prescribed exponent were studied by Hall--Higman \cite{HallHigman}, Kostrikin \cite{Kostrikin}, 
and Zel'manov \cite{Zelmanov1990,Zelmanov1991} in the context of the restricted 
Burnside problem. \\[0.1cm]
\item[(b)]
A group $G$ is $n$-abelian (in the sense of Baer \cite{Baer}) if and only if the linear 
polynomial $f(x):=-n+x$ is an ordered identity of some endomorphism $\gamma$ of $G$. These 
$n$-abelian groups were classified by Alperin \cite{Alperin}. \\[0.1cm]
\item[(c)]
Let $\alpha$ be an automorphism of a finite group $G$. If the polynomial $f(x):=-1+x^n$ is an 
ordered identity of $\alpha$, then $|\alpha|$ divides $n$. And vice-versa. 
Finite groups with a fixed-point-free automorphism of prescribed order have been studied 
extensively in the literature. Apart from the already mentioned results
\cite{Higman,KrekninKostrikin,Thompson1959a}, we highlight the work of 
Thompson \cite{Thompson1964}, Dade \cite{Dade}, Berger \cite{Berger}, and 
Jabara \cite{Jabara}. \\[0.1cm]
\item[(d)]
Let $\alpha$ be an automorphism of a finite group $G$. If the polynomial $f(x):=(x^n-
1)/(x-1)=1+x+\cdots+x^{n-1}$ is an ordered identity of $\alpha$, then we say that $\alpha$ is an $n$-splitting
automorphism. Such automorphisms have appeared in the work of Ersoy \cite{Ersoy}, Hughes--Thompson
\cite{HughesThompson} and Kegel \cite{Kegel}, Khukhro \cite{Khukhro1986}, and Zel'manov \cite{Zelmanov1992}.
We also mention Zel'manov's recent generalization \cite{Zelmanov2017}.
\end{itemize}
\end{ex}

By ignoring the order of the factors in Definition~\ref{1.1.1}, we obtain the more general notion of
identities of endomorphisms. These identities naturally appear as the characteristic polynomials of
endomorphisms of polycyclic groups.

\begin{defi}\label{1.1.3}
Let $G$ be a group with endomorphism $\gamma$. The polynomial $f(x)\in \Z[x]$ is an \emph{identity} of the
endomorphism if $f(x)$ can be decomposed into a finite sum $f(x)=f_1(x)+\cdots+f_k(x)$ of polynomials
$f_1(x),\ldots,f_k(x)\in \Z[x]$ such that the map $f_1(\gamma)\cdots f_k(\gamma):G\longrightarrow
G:g\longmapsto(f_1(\gamma)g)\cdots(f_k(\gamma)g)$ vanishes.
\end{defi}

\subsection{Main results}\label{sec:1.2}

\begin{defi}\label{1.2.1}
A polynomial $f(x)\in \Z[x]$ is \emph{powerful} if $f(x)\in \Z[x^2]\cup \Z[x^3]\cup \Z[x^4]\cup \cdots $.
\end{defi}

We can now formulate our first main result.

\begin{thm}[{Nilpotency}]\label{1.2.2}
For every polynomial $f(x)\in \Z[x]$ without powerful divisors of positive degree, there is a non-zero,
integer invariant $\iota_1(f(x))$ with the following property. Let $G$ be a finite group admitting a
fixed-point-free automorphism $\alpha$. Suppose that $f(x)$ is an identity of $\alpha$ and that
$\gcd(|G|,\iota_1(f(x)))=1$. Then $G$ is nilpotent.
\end{thm}

We will prove the theorem in several steps. We first define the invariant $\iota_1(f(x))$ (in
Definition~\ref{3.2.1}), we show that it is a non-zero integer (in Proposition~\ref{3.2.2}), and we show that
it describes the periodic decompositions of $f(x)$ (in Remark~\ref{3.1.2}). We then use the latter property to
generalize techniques of Higman \cite{Higman} in the context of the Frobenius conjecture. We will also use
Thompson's $p$-complement theorem of \cite{Thompson1959b}.

\begin{defi}\label{1.2.3}
A subset $X$ of an abelian group $(A,\cdot)$ is \emph{arithmetically-free} if $X$ contains no subset of the
form $\{a,a\cdot b,a\cdot b^2,a\cdot b^3,\ldots \}\cup \{b\}$.
\end{defi}

These sets were introduced in \cite{Moens} to generalize a classic result of Higman on the bounded nilpotency
of group-graded Lie rings. We can now formulate our second main result.

\begin{thm}[{Bounded class}]\label{1.2.4}
For every polynomial $f(x)\in \Z[x]$ for which the roots form an arithmetically-free subset $X$ of
$(\overline{\Q}^{\times},\cdot)$, there is a non-zero, integer invariant $\iota_2(f(x))$ with the following
property. Let $G$ be a nilpotent group and let $\alpha$ be an automorphism of $G$. Suppose that $f(x)$ is an
identity of $\alpha$ and that every finite subgroup $H$ of $G$ satisfies $\gcd(|H|,\iota_2(f(x)))=1$. Then the
nilpotency class of $G$ is at most $|X|^{2^{|X|}}$.
\end{thm}

We will again divide the proof into several steps. We first define the invariant $\iota_2(f(x))$ (in
Definition~\ref{3.2.1}), we show that it is a non-zero integer (in Proposition~\ref{3.2.2}), and we show that
it describes products of the roots of $f(x)$ (in Remark~\ref{3.1.6}). We then use the latter property to prove
a technical embedding theorem for Lie rings (Theorem~\ref{6.1.2}) by generalizing techniques of
Kreknin--Kostrikin \cite{KrekninKostrikin} and Higman \cite{Higman}. We next prove a grading theorem for Lie
rings (Theorem~\ref{6.2.1}), which generalizes theorems of Kreknin--Kostrikin \cite{KrekninKostrikin} and
Khukhro \cite{Khukhro2008} and which refines the main result of \cite{Moens}. We then derive a Lie ring
analogue of Theorem~\ref{1.2.4} (Theorem~\ref{6.3.1}) and Theorem~\ref{1.2.4} itself.

The assumptions of the two theorems are connected in a straight-forward manner.

\begin{prop}\label{1.2.5}
Let $f(x)\in \Z[x]$ with $f(0)\cdot f(1)\ne 0$. If $f(x)$ has no powerful divisors of positive degree, then
the roots of $f(x)$ form an arithmetically-free subset of $(\overline{\Q}^{\times},\cdot)$.
\end{prop}

It will then not be difficult to derive Corollary~\ref{1.0.1} from these two theorems. Let us mention one more
application of these results.

\begin{defi}\label{1.2.6}
Let $n\in \Z_{\geqslant 1}$. An automorphism $\alpha$ of a group $G$ is $n$\emph{-cyclotomic} if the
cyclotomic polynomial $\Phi_n(x)$ is an ordered identity of $\alpha$.
\end{defi}

Let $\varphi$ be the Euler totient function, let $\Kh $ be the function of Khukhro's theorem $8.2.1$ in
\cite{Khukhro1993}, and let $n_0$ be the square-free radical of a given natural number $n$.

\begin{cor}\label{1.2.7}
Let $n\in \Z_{\geqslant 2}$ and $k\in \Z_{\geqslant 1}$. Let $G$ be a finite group on $k$ generators and
suppose that $G$ has an $n$-cyclotomic automorphism. Then $G$ is nilpotent of class at most
$\varphi(n_0)^{2^{\varphi(n_0)}}$ if $n_0$ is composite, and $G$ is nilpotent of class at most $\max
\{\varphi(n_0)^{2^{\varphi(n_0)}},\Kh(k,n_0)\}$ if $n_0$ is a prime.
\end{cor}

This theorem implies, in particular, Kostrikin's solution \cite{Kostrikin} of the restricted Burnside problem
in prime exponent. Indeed, we need only specialize $n$ to a prime and $\alpha$ to the automorphism $x\mapsto
x$. We note, however, that the proof of the corollary implicitly relies on \cite{Kostrikin}, via subsequent
work of Khukhro \cite{Khukhro1986}, and on the work of Hughes--Thompson \cite{HughesThompson} and Kegel
\cite{Kegel}. Examples show, moreover, that the dependence of the upper bound on the number of generators
cannot be removed.

\subsection{Remarks}\label{sec:1.3}

Let us now justify the corollaries and their assumptions. We had already observed that the assumptions are
necessary. But the assumptions also suffice to recover a number of well-known results in the literature.
Indeed, Corollaries~\ref{1.0.1} and \ref{1.2.7} naturally extend the solution of the restricted Burnside
problem in prime exponent by Kostrikin (cf. Example~\ref{1.1.2} (a)), the classification of the finite
$n$-abelian groups by Baer and Alperin (cf. Example~\ref{1.1.2} (b)), the solutions of the Frobenius
conjecture by Thompson, Higman, and Kreknin--Kostrikin (cf. Example~\ref{1.1.2} (c)), and the work on
splitting automorphisms by Hughes--Thompson, Kegel, and Khukhro (cf. Example~\ref{1.1.2} (d)).

These classic results correspond with certain families of irreducible polynomials. More generally, every
monic, irreducible polynomial $f(x)\in \Z[x]\setminus \{x,x-1\}$ is an ordered identity of a fixed-point-free
automorphism of a finite, non-trivial group that is coprime to $\iota(f(x))$. So every such $f(x)$ appears in
Corollary~\ref{1.0.1}. Conversely, every fixed-point-free automorphism of a finite, non-trivial group has an
ordered identity that is monic and irreducible. We will prove precise and quantitative versions of these
statements in Section~\ref{sec:4}.

We conclude, in particular, that Corollary~\ref{1.0.1} can be applied to \emph{almost-every} monic polynomial
in the sense of probabilistic Galois theory. Indeed, almost-every monic polynomial is irreducible in the
bounded-degree model for $\Z[x]$. This classic result goes back to Hilbert and van der Waerden
\cite{VanDerWaerden}. Recent breakthroughs by Breuillard--Varj\'u \cite{BV}, Bary-Soroker--Kozma \cite{BK}, and
Bary-Soroker--Koukoulopoulos--Kozma \cite{BKK} show that, also in the bounded-height model for $\Z[x]$,
almost-every monic polynomial is irreducible. The same is true for various other models, such as the matrix
models of Rivin \cite{Rivin} and Eberhard \cite{Eberhard}.

\subsection{Overview}\label{sec:1.4}

In Section~\ref{sec:2}, we collect some results from the literature. In Section~\ref{sec:3}, we introduce and
study the \emph{invariants} $\iota_1(f(x))$ and $\iota_2(f(x))$. We show, in particular, how to compute the
invariants (Remark~\ref{3.1.8}) and we prove Proposition~\ref{1.2.5}. In Section~\ref{sec:4}, we establish the
\emph{existence} of identities by proving Proposition~\ref{4.2.1}, Proposition~\ref{4.2.2}, Lemma~\ref{4.2.4},
and Proposition~\ref{4.2.6}. We then turn to our \emph{structure} results. In Section~\ref{sec:5}, we prove
Theorem~\ref{1.2.2}. In section~\ref{sec:6}, we prove the auxiliary Theorems~\ref{6.1.2}, \ref{6.2.1}, and
\ref{6.3.1}, which may be of independent interest. We then prove Theorem~\ref{1.2.4} in Section~\ref{sec:7}.
In Section~\ref{sec:8}, we derive strong versions of Corollary~\ref{1.0.1} and Corollary~\ref{1.2.7}.

\section{Preliminaries}\label{sec:2}

\begin{lem}[{Gauss}]\label{2.0.1}
Let $f(x)\in \Z[x]$ and $g(x),h(x)\in \Q[x]$ with $f(x)=g(x)\cdot h(x)$. Then there exists some $a\in \Q
\setminus \{0\}$ such that $a\cdot g(x),a^{-1}\cdot h(x)\in \Z[x]$ and $f(x)= (a\cdot g(x))\cdot(a^{-1}\cdot
h(x))$.
\end{lem}

\begin{thm}[{Apostol \cite[Theorems 3--4]{Apostol}}]\label{2.0.2}
Let $m,n\in \Z$ with $1\leqslant m< n$. If $n/m$ is a natural power of a prime $q$, then
$\Res(\Phi_m(x),\Phi_n(x))=q^{\varphi(m)}$. Else, $\Res(\Phi_m(x),\Phi_n(x))=1$.
\end{thm}

\begin{defi}\label{2.0.3}
Let $n\in \Z \setminus \{0\}$ and let $G$ be a group. We say that $G$ is a $\pi(n)$\emph{-group} if the order
of every element of $G$ divides a natural power of $n$. A $\pi(n)$\emph{-subgroup} of $G$ is a subgroup of $G$
that is a $\pi(n)$-group. We say that $G$ is a $\pi(n)'$\emph{-group} if the only $\pi(n)$-subgroup of $G$ is
the trivial group.
\end{defi}

\begin{rem}\label{2.0.4}
A finite group $G$ is a $\pi(n)'$-group if and only if $\gcd(|G|,n)=1$.
\end{rem}

\begin{lem}\label{2.0.5}
Let $G$ be a finite group with a fixed-point-free automorphism $\alpha$.
\begin{itemize}
\item[(a)]
If $B$ is a normal, $\alpha$-invariant subgroup of an $\alpha$-invariant subgroup $A$ of $G$, then the induced
automorphism $\alpha_{A/B}:A/B\longrightarrow A/B:x\cdot B\longmapsto \alpha(x)\cdot B$ is also
fixed-point-free.
\item[(b)]
For every prime $p$, there is a Sylow $p$-subgroup $P$ of $G$ satisfying $\alpha(P)=P$.
\item[(c)]
Suppose that $G$ is solvable. Then, for every $n\in \Z \setminus \{0\}$, there is a Hall $\pi(n)$-subgroup $H$
of $G$ satisfying $\alpha(H)=H$.
\end{itemize}
\end{lem}

\begin{proof}
(a) This can be found as theorem $1.6.1$ of \cite{Khukhro1993}. (b) According to Sylow's second theorem, there
is a $g\in G$ such that $\alpha(P)=P^g$. Since $\alpha$ is injective, the map $h\mapsto h\cdot \alpha(h^{-1})$
is surjective. So we can select $h\in G$ with $\alpha(g)=h\cdot \alpha(h^{-1})$. Then $\alpha(P^h)=P^h$. (c)
We need only replace Sylow's theorem with Hall's main theorem of \cite{Hall} in the previous proof.
\end{proof}

\begin{thm}[{Rowley \cite{Rowley}}]\label{2.0.6}
If a finite group $G$ has a fixed-point-free automorphism, then $G$ is solvable.
\end{thm}

For a subgroup $H$ of a group $G$ we let $\mathrm{N}_G(H)$ resp. $\mathrm{C}_G(H)$ be the normalizer resp.
centralizer of $H$ in $G$.

\begin{thm}[{Thompson \cite[Theorem A]{Thompson1959b}}]\label{2.0.7}
Let $p$ be an odd prime and let $G$ be a finite group with a group $A$ of automorphisms. Let $P$ be a Sylow
$p$-subgroup of $G$ and suppose that $P$ is invariant under $A$. Suppose that, for every normal, $A$-invariant
subgroup $Q$ of $P$, all the elements of $\mathrm{N}_G(Q)$ that have order coprime to $p$ are in
$\mathrm{C}_G(Q)$. Then $G$ has a characteristic subgroup $K$ such that $G\cong K\rtimes P$.
\end{thm}

\begin{thm}[{Gruenberg \cite[Theorem 2.1]{Gruenberg}}]\label{2.0.8}
Let $G$ be a finitely-generated, torsion-free, nilpotent group. Then, for every prime $p$, the group is
residually a finite $p$-group.
\end{thm}

\begin{ex}[{Moens \cite[Example 2.2]{Moens}}]\label{2.0.9}
If $X$ is the root set of an irreducible polynomial $f(x)\in \Z[x]\setminus \{\pm x,\pm(x-1)\}$, then $X$ is
an arithmetically-free subset of $(\overline{\Q}^{\times},\cdot)$.
\end{ex}

\begin{thm}[{Moens \cite[Theorem 3.14]{Moens}}]\label{2.0.10}
For every finite, arithmetically-free subset $X$ of an arbitrary group $(F,\cdot)$, there exists a minimal
integer $H(X,F)$ with the following property. If a Lie ring $L$ is graded by $F$ and supported by $X$, then
$L$ is nilpotent of class at most $H(X,F)$.
\end{thm}

In \cite{Moens}, the theorem is stated for Lie algebras with coefficients in a field, but the proofs are valid
for Lie rings with coefficients in $\Z$.

\begin{thm}[{Kreknin--Kostrikin \cite{KrekninKostrikin}; Shalev \cite[Proposition 2.4]{Shalev}}]\label{2.0.11}
Let $X$ be a finite subset of the multiplicative group $(\F^{\times},\cdot)$ of a field $\F$ and suppose that
$1\notin X$. Then every Lie ring $L$ that is graded by $\F^{\times}$ and supported by $X$ is solvable of
derived length at most $2^{|X|}$.
\end{thm}

\begin{thm}[{Hughes--Thompson \cite{HughesThompson}; Kegel \cite[Satz 1]{Kegel}}]\label{2.0.12}
Let $p$ be a prime and let $G$ be a finite group. If $G$ has a $p$-cyclotomic automorphism, then $G$ is
nilpotent.
\end{thm}

\begin{thm}[{Khukhro \cite{Khukhro1986}; \cite[Theorem 8.2.1]{Khukhro1993}}]\label{2.0.13}
Let $p$ be a prime and let $k\in \Z_{\geqslant 1}$. Then there is a natural number $\Kh(k,p)$ with the
following property. Every finite $p$-group on $k$ generators admitting a $p$-cyclotomic automorphism is
nilpotent of class at most $\Kh(k,p)$.
\end{thm}

\begin{thm}[{Peter--Weyl \cite{PeterWeyl}; \cite[Theorem 1.4.12]{Tao}}]\label{2.0.14}
Every compact Hausdorff group is residually a linear group over a field of characteristic $0$.
\end{thm}

\begin{thm}[{Mal'cev \cite[Theorem VI]{Malcev1965}}]\label{2.0.15}
Every finitely-generated, linear group over a field of characteristic $0$ is residually a finite group.
\end{thm}

\section{The invariants and their properties}\label{sec:3}

\subsection{\texorpdfstring{The auxiliary invariants $\Cong(f(x))$, $\Prod(f(x))$, and $\Discr_{\ast}(f(x))$}
  {The auxiliary invariants Cong(f(x)), Prod(f(x)), and Discr*(f(x))}}\label{sec:3.1}

\begin{defi}\label{3.1.1}
Let $f(x)=a_0+a_1\cdot x+\cdots+a_d\cdot x^d\in \Z[x]$ and let $n\in \Z_{\geqslant 1}$. For each $j\in
\{0,1,\ldots,n-1\}$, we define the partial sums
\[
f_{n,j}(x):=\sum_{i\equiv j\bmod n}a_i\cdot x^i,
\]
so that $f_{n,j}(x)\in x^j\cdot \Z[x^n]$ and $f(x)=f_{n,0}(x)+\cdots+f_{n,n-1}(x)$. We next define the ideal
\[
I_n:=f_{n,0}(x)\cdot \Z[x]+\cdots+f_{n,n-1}(x)\cdot \Z[x]
\]
of $\Z[x]$ that is generated by the $f_{n,j}(x)$. We then define the principal ideal $I:=\Z \cap I_2\cap
I_3\cap \cdots=(\Z \cap I_2)\cap(\Z \cap I_3)\cap \cdots $ of $\Z$ and we let $\Cong(f(x))$ be the
non-negative generator of this $I$.
\end{defi}

\begin{rem}\label{3.1.2}
For each $n\geqslant 2$, there exist $g_0(x),\ldots,g_{n-1}(x)\in \Z[x]$ such that
\[
\Cong(f(x))=g_0(x)\cdot f_{n,0}(x)+\cdots+g_{n-1}(x)\cdot f_{n,n-1}(x).
\]
\end{rem}

Let us now characterize the polynomials $f(x)$ for which $\Cong(f(x))\ne 0$.

\begin{lem}\label{3.1.3}
Let $a(x),b(x)\in \Q[x]$ and suppose that they are not both $0$. Let $n\in \Z_{\geqslant 1}$ and define
$h(x):=\gcd_{\Q[x]}(a(x),b(x))$ and $H(x):=\gcd_{\Q[x]}(a(x^n),b(x^n))$. Then $h(x^n)=H(x)$.
\end{lem}

\begin{proof}
We abbreviate $A(x):=a(x^n)$ and $B(x):=b(x^n)$. By definition, there exists $f(x),g(x)\in \Q[x]$ such that
$a(x)=h(x)\cdot f(x)$ and $b(x)=h(x)\cdot g(x)$. By substituting $x\mapsto x^n$, we obtain $A(x)=h(x^n)\cdot
f(x^n)$ and $B(x)=h(x^n)\cdot g(x^n)$. So $h(x^n)|H(x)$ in $\Q[x]$. Moreover, according to Bezout's theorem,
there exist $v(x),w(x)\in \Q[x]$ such that $h(x)=v(x)\cdot a(x)+w(x)\cdot b(x)$. By substitution, we obtain
$h(x^n)=v(x^n)\cdot A(x)+ w(x^n)\cdot B(x)$. So $H(x)|h(x^n)$ in $\Q[x]$. Since $h(x^n)$ and $H(x)$ are monic
by definition, we conclude $h(x^n)=H(x)$.
\end{proof}

\begin{prop}\label{3.1.4}
Let $f(x)\in \Z[x]$ with $f(0)\ne 0$. Then the following are equivalent:
\begin{itemize}
\item[(a)]
$\Cong(f(x))\ne 0$.
\item[(b)]
All powerful divisors of $f(x)$ are constant.
\end{itemize}
\end{prop}

\begin{proof}
Suppose first that $\Cong(f(x))\ne 0$. Let $g(x^n)$ be a powerful divisor of $f(x)$, for some $g(x)\in \Z[x]$
and $n\geqslant 2$. By definition, there exists some $h(x)\in \Z[x]$ such that $f(x)= h(x)\cdot g(x^n)$. So we
have $f(x)=h_{n,0}(x)\cdot g(x^n)+\cdots+h_{n,n-1}(x)\cdot g(x^n)$. Since $h_{n,j}(x)\cdot g(x^n)$ belongs to
$x^j\cdot \Z[x^n]$, we may conclude that $f_{n,j}(x)=h_{n,j}(x)\cdot g(x^n)$. According to Remark~\ref{3.1.2},
there exist $h_0(x),\ldots,h_{n-1}(x)\in \Z[x]$ such that $h_0(x)\cdot f_{n,0}(x)+\cdots+ h_{n-1}(x)\cdots
f_{n,n-1}(x)=\Cong(f(x))$. So $g(x^n)|\Cong(f(x))$ and the powerful divisor $g(x^n)$ of $f(x)$ is constant.

Suppose next that $f(x)$ has no powerful divisors of positive degree and set $d:= \deg(f(x))$. Since
$f_{n,j}(x)=0$ for all $j\geqslant d+1$, we have $I_{d+1}=I_{d+2}=\cdots $. So we need only show that the
$I_n$ are non-zero for $n\geqslant 2$. So let $n\in \Z_{\geqslant 2}$. For each $0\leqslant j\leqslant n-1$,
there is a unique $F_{n,j}(x)\in \Z[x]$ such that $f_{n,j}(x):=F_{n,j}(x^n)\cdot x^j$. Since $x$ does not
divide $f(x)$, the greatest common divisor $H(x)$ of the $f_{n,j}(x)$ in $\Q[x]$ coincides with the greatest
common divisor of the $F_{n,j}(x^n)$ in $\Q[x]$. According to Lemma~\ref{3.1.3}, there then exists some
$h(x)\in \Q[x]$ such that $H(x)=h(x^n)$. According to Gauss' lemma~\ref{2.0.1}, some scalar multiple of
$h(x^n)$ is in $\Z[x^n]$ and divides $f(x)$ in $\Z[x]$. By assumption, the latter can only occur if $h(x^n)$
is constant. So $h(x^n)\in \Z \setminus \{0\}$ and the polynomials $f_{n,j}(x)$ are coprime over $\Q$.
According to Bezout's theorem, we then have $\Z \cap I_n\ne \{0\}$.
\end{proof}

\begin{defi}\label{3.1.5}
Let $f(x)\in \Z[x]$. If $f(x)$ is constant, then we define $\Discr_{\ast}(f(x)):= f(x)$ and
$\Prod(f(x)):=f(x)^2$. Else, we proceed as follows. We consider the factorization
\[
f(x)=a\cdot \prod_{1\leqslant i\leqslant l}(x-\lambda_i)^{m_i},
\]
of $f(x)$, where $a$ is the leading coefficient, and where $\lambda_1,\ldots,\lambda_l$ are the distinct roots
with respective multiplicities $m_1,\ldots,m_l$. We define $m:=\max \{m_1,\ldots,m_l\}$ and $d:=
\deg(f(x))=m_1+\cdots+m_l$. We then define
\[
\begin{aligned}
\Discr_{\ast}(f(x))
  &:= a^{1+2d^2}\cdot(m-1)!\cdot
      \prod_{\substack{1\leqslant i,j\leqslant l\\ i\ne j}}
      (\lambda_i-\lambda_j)^m,\quad\text{and}\\
\Prod(f(x))
  &:= a^{2d^3}\cdot
      \prod_{\substack{1\leqslant i,j\leqslant l\\ f(\lambda_i\cdot\lambda_j)\ne0}}
      f(\lambda_i\cdot\lambda_j)\\
  &= a^{2d^3}\cdot
      \prod_{\substack{1\leqslant i,j,k\leqslant l\\ f(\lambda_i\cdot\lambda_j)\ne0}}
      a\cdot(\lambda_i\cdot\lambda_j-\lambda_k)^{m_k}.
\end{aligned}
\]
\end{defi}

\begin{rem}\label{3.1.6}
Suppose $f(x)\notin \Z$ and let $i,j\in \{1,\ldots,l\}$. If $\lambda_i\ne \lambda_j$, then there is some $u\in
\Z[\lambda_1,\ldots,\lambda_l]$ and some $k\in \Z_{\geqslant 1}$ such that
$u\cdot(\lambda_i-\lambda_j)=\Discr_{\ast}(f(x))^k$. If $f(\lambda_i\cdot \lambda_j)\ne 0$, then there is some
$u\in \Z[\lambda_1,\ldots,\lambda_l]$ and some $k\in \Z_{\geqslant 1}$ such that $u\cdot f(\lambda_i\cdot
\lambda_j)= \Prod(f(x))^k$.
\end{rem}

\begin{prop}\label{3.1.7}
For all $f(x)\in \Z[x]\setminus \{0\}$, we have $\Discr_{\ast}(f(x)),\Prod(f(x))\in \Z \setminus \{0\}$.
\end{prop}

\begin{proof}
We continue to use the notation of the definition. Since $f(x)\ne 0$, the invariants are non-zero by
definition. So we need only verify that they are integers. We may assume that $d\geqslant 2$, since the cases
$d\leqslant 2$ are straightforward to check. Let $u(x)$ be the largest square-free factor of $f(x)$ in
$\Z[x]$, say with leading factor $\overline{a}\in \Z$. Define $v(x):=(a/\overline{a})\cdot u(x)$. Then
$v(x)=a\cdot(x-\lambda_1)\cdots(x-\lambda_l)$. Since $\overline{a}$ divides $a$, we also have $v(x)\in \Z[x]$.
We see, in particular, that $\Discr(v(x))\in \Z[x]$. Now set $k:=1+2d^2-2m(l-1)\in \Z_{\geqslant 0}$. Then it
is not difficult to verify that $\Discr_{\ast}(f(x))=\pm(m-1)!\cdot a^k\cdot \Discr(v(x))^m$. So we may
conclude that $\Discr_{\ast}(f(x))\in \Z$. Now define the auxiliary polynomials
\[
\begin{aligned}
\widetilde g(x)&:=\prod_{\substack{1\leqslant i,j\leqslant l\\ f(\lambda_i\cdot\lambda_j)=0}}
  (x-\lambda_i\cdot\lambda_j) \quad\text{and}\\
\widetilde h(x)&:=\prod_{\substack{1\leqslant i,j\leqslant l\\ f(\lambda_i\cdot\lambda_j)\ne0}}
  (x-\lambda_i\cdot\lambda_j).
\end{aligned}
\]
Let $E$ be the splitting field over $\Q$ of the separable polynomial $v(x)$. Then every element in
$\Gal(E/\Q)$ induces a permutation of the roots of $\widetilde{g}(x)$. Since the coefficients of
$\widetilde{g}(x)$ are elementary symmetric in the roots of $\widetilde{g}(x)$, we conclude that
$\widetilde{g}(x)\in \Q[x]$. Similarly, we obtain that $\widetilde{h}(x)\in \Q[x]$. Now note that
$\widetilde{g}(x)\cdot \widetilde{h}(x)=\prod_{1\leqslant i,j\leqslant l}(x-\lambda_i\cdot \lambda_j)=
\chi_{C\otimes C}(x)$, where $C$ is the companion operator of $(1/a)\cdot v(x)$. Since $a\cdot C\in
\End(\Z^d)$, we have $a^{2l^2}\cdot \widetilde{g}(x)\cdot \widetilde{h}(x)\in \Z[x]$. According to Gauss'
lemma~\ref{2.0.1}, we then have $h(x):= a^{2l^2}\cdot \widetilde{h}(x)\in \Z[x]$. So $\Res(h(x),f(x))\in \Z$.
Now set $k:=(2d^3-2l^2d)$. Then it is not difficult to verify that $\Prod(f(x))=\pm a^k\cdot \Res(h(x),f(x))$.
We may therefore conclude that also $\Prod(f(x))\in \Z$.
\end{proof}

\begin{rem}\label{3.1.8}
The proof shows us how to compute the invariants by means of resultants and without having to know the roots
of $f(x)$. Moreover, $\Discr_{\ast}(f(x))$ clearly coincides with $\Discr(f(x))$ (up to a sign) if $f(x)$ is
monic and square-free.
\end{rem}

\subsection{\texorpdfstring{The invariants $\iota_1$ and $\iota_2$ of $f(x)$}
  {The invariants iota 1 and iota 2 of f(x)}}\label{sec:3.2}

\begin{defi}\label{3.2.1}
Let $f(x)\in \Z[x]\setminus \{0\}$. Consider the unique $\overline{f}(x)\in \Z[x]\setminus \{0\}$ and
$m_0,m_1\in \Z_{\geqslant 0}$ such that
\[
f(x)=x^{m_0}\cdot(x-1)^{m_1}\cdot \overline{f}(x)
\]
and such that $\overline{f}(0)\cdot \overline{f}(1)\ne 0$. We then define
\[
\begin{aligned}
\iota_1&:=\iota_1(f(x)):=\Cong(\overline f(x)),\\
\iota_2&:=\iota_2(f(x)):=\Discr_\ast(\overline f(x))\cdot\Prod(\overline f(x)).
\end{aligned}
\]
\end{defi}

\begin{prop}\label{3.2.2}
Let $f(x)\in \Z[x]$. (a) If $\overline{f}(x)$ has no powerful divisors of positive degree, then
$\iota_1(f(x))\in \Z \setminus \{0\}$. (b) If $f(x)\ne 0$, then $\iota_2(f(x))\in \Z \setminus \{0\}$.
\end{prop}

\begin{proof}
(a) As a divisor of the non-zero polynomial $f(x)$, the non-zero polynomial $\overline{f}(x)$ has no powerful
divisors of positive degree. Proposition~\ref{3.1.4} then implies that $\Cong(\overline{f}(x))\ne 0$. (b) We
need only apply Proposition~\ref{3.1.7}.
\end{proof}

\subsection{\texorpdfstring{The invariants $\iota_1$ and $\iota_2$ of $\Phi_n(x)$}
  {The invariants iota 1 and iota 2 of Phi n(x)}}\label{sec:3.3}

In order to compute the prime divisors of $\iota_1(\Phi_n(x))$ and $\iota_2(\Phi_n(x))$, we recall some
well-known properties of cyclotomic polynomials. Let $n\in \Z_{>1}$ and let $n_0$ be the \emph{square-free
radical} of $n$, i.e.: the largest (positive) square-free factor of $n$. Let $p$ be an odd prime not dividing
$n_0$. Then $\deg(\Phi_n(x))=\varphi(n)$, where $\varphi$ is the Euler totient-function. Moreover,
\begin{equation}\label{eq:3.1}
\Phi_n(x)=\Phi_m(x^{n/m})
\end{equation}
and
\begin{equation}\label{eq:3.2}
\Phi_{p\cdot n_0}(x)\cdot \Phi_{n_0}(x)=\Phi_{n_0}(x^p).
\end{equation}
If $n_0$ is odd, then
\begin{equation}\label{eq:3.3}
\Phi_{2\cdot n_0}(x)=\Phi_{n_0}(-x).
\end{equation}
If $n_0$ is a prime, then $\Phi_n(1)=n_0$. Else, $\Phi_n(1)=1$. We continue to use the notation of
Definition~\ref{3.1.1}.

\begin{defi}\label{3.3.1}
Let $f(x)\in \Z[x]\setminus \{0\}$ and let $n\in \Z_{\geqslant 2}$. We define $\RRes_n(f(x))$ to be the
(unique) non-negative generator of the (principal) ideal $\Z \cap I_n$ of $\Z$.
\end{defi}

So $\Cong(f(x))$ is the least common multiple of $\RRes_2(f(x)),\ldots,\RRes_{\deg(f(x))+2}(f(x))$. We will
need two elementary properties of these reduced resultants.

\begin{lem}[{Division}]\label{3.3.2}
Let $a(x),b(x),c(x)\in \Z[x]$ and suppose that $a(x)\cdot b(x)= c(x)$. For every $n\in \Z_{\geqslant 2}$, we
have $\RRes_u(a(x))|\RRes_u(c(x))$.
\end{lem}

\begin{proof}
We note that for every $i\in \Z_{\geqslant 0}$, we have $c_{n,i}(x)=\sum_{j+k\equiv i\bmod n}a_{n,j}(x)\cdot
b_{n,k}(x)$. By definition, there exist $C_0(x),\ldots,C_{u-1}(x)\in \Z[x]$ such that $\RRes_u(c(x))=
\sum_{0\leqslant i<n}C_i(x)\cdot c_{n,i}(x)$. So
\[
\RRes_n(c(x))=\sum_{0\leqslant i<n}C_i(x)\cdot
  \left(\sum_{j+k\equiv i\bmod n}a_{n,j}(x)\cdot b_{n,k}(x)\right)
\]
and therefore $\RRes_n(c(x))\in(a_{n,0}(x)\cdot \Z[x]+\cdots+a_{n,n-1}(x)\cdot \Z[x])\cap \Z$.
\end{proof}

\begin{lem}[{Composition}]\label{3.3.3}
Consider $f(x)\in \Z[x]\setminus \{0\}$. Consider integers $n\geqslant 2$, and $k,l\geqslant 1$. If $l\cdot
k\equiv 1\bmod n$, then $\RRes_n(f(x^l))|\RRes_n(f(x)).$
\end{lem}

\begin{proof}
Let the polynomial be given by $f(x)=\sum_{0\leqslant j\leqslant d}a_j\cdot t^j$ and define $s(x):=f(x^l)$.
Then for all $0\leqslant i<n$, we have
\[
s_{n,i}(x)=\sum_{\substack{0\leqslant j\leqslant d\\ l\cdot j\equiv i\bmod n}}
 a_j\cdot(x^l)^j
 =\sum_{\substack{0\leqslant j\leqslant d\\ j\equiv k\cdot i\bmod n}}a_j\cdot(x^l)^j
 =f_{n,k\cdot i\bmod n}(x^l).
\]
By definition, there exist $c_{n,0}(x),\ldots,c_{n,n-1}(x)\in \Z[x]$ such that $\RRes_n(f(x))=
\sum_jc_{n,j}(x)\cdot f_{n,j}(x)$. By substituting $x\mapsto x^l$, we obtain
$\RRes_n(f(x))=\sum_jc_{n,j}(x^l)\cdot f_{n,j}(x^l)\in s_{n,0}(x)\cdot \Z[x]+\cdots+s_{n,n-1}(x)\cdot \Z[x].$
So $\RRes_n(r(x))\in(s_{n,0}(x)\cdot \Z[x]+\cdots+ s_{n,n-1}(x)\cdot \Z[x])\cap \Z$.
\end{proof}

Let us now fix $n\in \Z_{\geqslant 2}$. The factorization $\Phi_n(x)=x^{m_0}\cdot(x-1)^{m_1}\cdot
\overline{\Phi}_n(x)$ of Definition~\ref{3.2.1} is rather trivial: $m_0=m_1=0$ and
$\overline{\Phi}_n(x)=\Phi_n(x)$. So $\iota_1(\Phi_n(x))=\Cong(\Phi_n(x))$ and
$\iota_2(\Phi_n(x))=\Discr_{\ast}(\Phi_n(x))\cdot \Prod(\Phi_n(x))$. Let $n_0$ be the square-free radical of
$n$.

\begin{prop}\label{3.3.4}
Let $n\in \Z_{\geqslant 2}$. If $n\ne n_0$, then $\iota_1(\Phi_n(x))=1$. Else, $\iota_1(\Phi_n(x))=0$.
\end{prop}

\begin{proof}
If $n\ne n_0$, then \eqref{eq:3.1} tells us that $\Phi_{n_0}(x^{n/n_0})$ is a powerful divisor of $\Phi_n(x)$
of positive degree. In this case, Proposition~\ref{3.1.4} tells us that $\Cong(\Phi_n(x))=0$. So we may from
now on assume that $n=n_0$. We then need only show that $\RRes_u(\Phi_n(x))=1$ for all $2\leqslant u\leqslant
\varphi(n)+1$. We proceed in three steps. $(i.)$ Suppose that $n$ is a prime. We may assume that $n\ne 2$,
since otherwise the proof is trivial. Set $f(x):=\Phi_n(x)$ in order to avoid confusing subscripts. Suppose
first that $u=\varphi(n)+1=(n-1)+ 1=n$. Then $f_{u,0}(x)=1$, so that $1=1\cdot f_{u,0}(x)+0\cdot
f_{u,1}(x)+\cdots+0\cdot f_{u,u-1}(x)$, and therefore $\RRes_u(r(x))=1$. We next consider $2\leqslant u<n$.
Then $n-1\not\equiv u-1\bmod u$, so that $1=1\cdot f_{u,0}(x)+0\cdot g_{u,1}(x)+\cdots+0\cdot
f_{u,u-2}(x)-x\cdot f_{u,u-1}(x).$ This also implies that $\RRes_u(f(x))=1$. We conclude that
$\Cong(\Phi_n(x))=1$. $(ii.)$ Suppose next that $n$ is odd. Let us proceed by induction on the number $l$ of
distinct prime factors of $n$. The base of the induction, $l=1$, is given by the previous paragraph. So we
suppose that $l>1$. Let $n=p_1\cdots p_l$ be the decomposition of $n$ into (distinct, odd) primes. We select
an integer $u$ satisfying $2\leqslant u\leqslant \varphi(n)+1=(p_1-1)\cdots(p_l-1)+1$. Then we note that there
is at least one $i\in \{1,\ldots,l\}$ such that $p_i$ does \emph{not} divide $u$. Formula \eqref{eq:3.2} tells
us that $\Phi_n(x)|\Phi_{n/p_i}(x^{p_i})$. Lemma~\ref{3.3.2}, lemma~\ref{3.3.3}, and the induction hypothesis
then imply that $\RRes_u(\Phi_n(x))|\RRes_u(\Phi_{n/p_i}(x^{p_i}))|\RRes_u(\Phi_{n/p_i}(x))=1.$ We conclude
that also in this case $\Cong(\Phi_n(x))=1$. $(iii.)$ Finally, we consider the even $n$. We may again assume
that $n\ne 2$, since otherwise the proof is trivial. Formula \eqref{eq:3.3} gives us the equality
$\Phi_n(x)=\Phi_{n/2}(-x)$. The odd case then tells us that $\RRes_u(\Phi_n(x))=\RRes_u(\Phi_{n/2}(-x))=
\RRes_u(\Phi_{n/2}(x))=1.$ We conclude $\Cong(\Phi_n(x))=1$.
\end{proof}

\begin{prop}\label{3.3.5}
The invariant $\iota_2(\Phi_n(x))$ divides a natural power of $n$.
\end{prop}

\begin{proof}
We assume that $n>2$, since otherwise the computation is straight-forward. Since the cyclotomic field $K_n$
corresponding with $\Phi_n(x)$ is monogenic, $\Discr(\Phi_n(x))$ coincides with the field discriminant of
$K_n$: $\Discr(\Phi_n(x))=\Delta_{K_n}= (-1)^{\varphi(n)/2}\cdot n^{\varphi(n)}/\prod
p^{\varphi(n)/\varphi(p)}$, where $p$ runs over the prime divisors of $n$. Since $\Phi_n(x)$ is monic and
square-free, we conclude that $\Discr_{\ast}(\Phi_n(x))=\pm \Discr(\Phi_n(x))$ divides $n^n$ (cf.
Remark~\ref{3.1.8}). Now consider $\Prod(\Phi_n(x))$. We use the polynomial $h(x)\in \Z[x]$ of the proof of
Proposition~\ref{3.1.7}. Each root of $h(x)$ is a \emph{non}-primitive $n$'th root of unity, so that $h(x)$
factorizes as $h(x)=\prod \Phi_m(x)^{a_m},$ for suitable $a_m\in \Z_{\geqslant 0}$, where $m$ runs over the
\emph{proper} divisors of $n$. We therefore have $\Prod(\Phi_n(x))=\Res(\Phi_n(x),h(x))= \prod
\Res(\Phi_n(x),\Phi_m(x))^{a_m}$, where $m$ runs over the \emph{proper} divisors of $n$. The factors are given
in Theorem~\ref{2.0.2}. So also $\Prod(\Phi_n(x))$ divides a natural power of $n$.
\end{proof}

\subsection{Arithmetically-free subsets of groups}\label{sec:3.4}

The polynomials of Theorem~\ref{1.2.2} have no powerful divisors of positive degree, while the polynomials of
Theorem~\ref{1.2.4} have an arithmetically-free set of roots. Let us show that, under mild conditions, the
former property implies the latter.

\begin{prop}\label{3.4.1}
Let $f(x)\in \Z[x]$ and suppose that $f(0)\ne 0$. Then the following are equivalent.
\begin{itemize}
\item[(a)]
The roots of $f(x)$ do not form an arithmetically-free subset of $(\overline{\Q}^{\times},\cdot)$.
\item[(b)]
There exists some $n\in \Z_{\geqslant 1}$ and $g(x)\in \Z[x]\setminus \Z$ such that $g(x^n)|f(x)$ and
$\Phi_n(x)|f(x)$.
\end{itemize}
\end{prop}

\begin{proof}
Let $X$ be the root set. Suppose first that $g(x^n),\Phi_n(x)|f(x)$. Let $\lambda$ be a root of $g(x^n)$ and
let $\mu$ be a root of $\Phi_n(x)$. Then $\lambda,\lambda \cdot \mu,\lambda \cdot \mu^2,\ldots $ are all roots
of $f(x)$. So, in this case, $X$ is not arithmetically-free.

Now suppose that $X$ is not arithmetically-free. Then there exist $\lambda,\mu \in X$ such that
$\lambda,\ldots,\lambda \cdot \mu^{|X|}\in X$. Then $\mu$ has finite order, say $n\geqslant 1$, and
$\Phi_n(x)|f(x)$ in $\Z[x]$. We now consider the partial decomposition $f(x)=f_{n,0}(x)+\cdots+f_{n,n-1}(x)$
of $f(x)$ and evaluate in the terms of the progression. We then obtain the linear Vandermonde system $\{\sum
\mu^{i\cdot j}\cdot f_{n,j}(\lambda)=0\}_{0\leqslant i,j\leqslant n-1}$. Since the order of $\mu$ is $n$, the
determinant of this system is non-zero. So the partial sums $f_{n,j}(\lambda)$ all vanish. Since $\lambda \ne
0$, we see that $\lambda$ is a common root of $f_{n,0}(x)\cdot x^{-0},f_{n,1}(x)\cdot
x^{-1},\ldots,f_{n,n-1}(x)\cdot x^{-(n-1)}$ and their greatest common divisor over $\Q$, say $h(x)$. Since
$f_{n,j}(x)\cdot x^{-j}\in \Z[x^n]$, Lemma~\ref{3.1.3} implies that there is a $\widetilde{g}(x)\in
\Q[x]\setminus \Q$ such that $h(x)=\widetilde{g}(x^n)$. According to Gauss' lemma~\ref{2.0.1}, a suitable
scalar multiple of $\widetilde{g}(x^n)$, say $g(x^n)$, is then a non-constant divisor of $f(x)$ in $\Z[x]$.
\end{proof}

\begin{prop}\label{3.4.2}
Let $f(x)\in \Z[x]$ with $f(0)\cdot f(1)\ne 0$. Suppose that $f(x)$ has no powerful divisors of positive
degree. Then the roots of $f(x)$ form an arithmetically-free subset of $(\overline{\Q}^{\times},\cdot)$.
\end{prop}

\begin{proof}
We suppose the statement is false and we will derive a contradiction. According to Proposition~\ref{3.4.1},
there is a $g(x)\in \Z[x]\setminus \Z$ and $n\geqslant 1$ such that $\Phi_n(x),g(x^n)|f(x)$. Since $f(x)$ has
no powerful divisors of positive degree, we have $n=1$, and therefore $f(1)=0$.
\end{proof}

\Needspace{7\baselineskip}
\section{Identities}\label{sec:4}

\subsection{General properties}\label{sec:4.1}

\begin{lem}\label{4.1.1}
Let $G$ be a group and let $\gamma$ be an endomorphism. Let $f(x),g(x)\in \Z[x]$ and define $h(x):=f(x)\cdot
g(x)$. Then there exist $k\in \Z_{\geqslant 0}$ and $h_1(x),\ldots,h_k(x)\in \Z[x]$, such that
$h(x)=h_1(x)+\cdots+h_k(x)$ and $(f(\gamma)\circ g(\gamma))(t)=(h_1(\gamma)(t))\cdots(h_k(\gamma)(t))$ for all
$t\in G$.
\end{lem}

\begin{proof}
Let $f(x):=a_0+a_1\cdot x+\cdots+a_d\cdot x^d$ and $g(x):=b_0+b_1\cdot x+\cdots+b_D\cdot x^D$. We proceed by
induction on the absolute height $l:=|a_0|+\cdots+|a_d|$ of $f(x)$. Suppose first that $l=0$. Then it suffices
to define $k=1$ and $h_1(x):=0$. If $l>0$, then we proceed as follows. Suppose first that $f(x)$ is a
monomial, say $a_j\cdot x^j$. If $a_j>0$, then we define $k:=|a_j|$ and $h_1(x):=\cdots:=h_k(x):=g(x)$. If
$a_j<0$, then we define $k:=|a_j|\cdot(D+1)$ and $h_{i+l\cdot D}(x):=-b_{D-i}\cdot t^{j+D-i}$, where $i\in
\{0,\ldots,D\}$ and $l\in \{0,\ldots,|a_j|-1\}$. So we may suppose that there is some $l$ such that
$F_1(x):=a_0+\cdots+ a_l\cdot x^l$ and $F_2(x):=a_{l+1}\cdot x^{l+1}+\cdots+a_d\cdot x^d$ have strictly
smaller absolute height than $f(x)$. The induction hypotesis then gives
$h_1(x),\ldots,h_{k_1}(x),h_{k_1+1}(x),\ldots,h_{k_1+k_2}(x)\in \Z[x]$ such that $\sum_{1\leqslant i\leqslant
k_1}h_i(x)=F_1(x)\cdot g(x)$ and $\sum_{1\leqslant i\leqslant k_2}h_{k_1+i}(x)=F_2(x)\cdot g(x)$ and such
that, for all $t\in G$, we have $(F_1(\gamma)\circ g(\gamma))(x)=(h_1(\gamma)(x))\cdots(h_{k_1}(\gamma)(x))$
and $(F_2(\gamma)\circ g(\gamma))(x)=(h_{k_1+1}(\gamma)(x))\cdots(h_{k_1+k_2}(\gamma)(x))$. So
$(f(\gamma)\circ g(\gamma))(x)=(F_1(\gamma)\circ g(\gamma)(x))\cdot(F_2(\gamma)\circ
g(\gamma)(x))=(h_1(\gamma)(x))\cdots(h_{k_1+k_2}(\gamma)(x))$ with $h_1(x)+\cdots+h_{k_1+k_2}(x)=
F_1(x)+F_2(x)=f(x)$.
\end{proof}

\begin{lem}\label{4.1.2}
Let $G$ be a group, let $\gamma$ be an endomorphism, and let $H$ be a normal subgroup of $G$ that is invariant
under $H$. Let $f_{G/H}(x)$ be an identity of the natural endomorphism $\gamma_{G/H}:G/H\longrightarrow
G/H:t\cdot H\longmapsto \gamma(t)\cdot H$ induced on $G/H$ and let $f_H(x)$ be an identity of the natural
endomorphism $\gamma_H:H\longrightarrow H$ on $H$, obtained by restriction. Then $f(x):=f_{G/H}(x)\cdot
f_H(x)$ is an identity of $\gamma$.
\end{lem}

\begin{proof}
The map $f_{G/H}(\gamma):G\longrightarrow G$ maps $G$ into $H$ and the map $f_H(\gamma):G\longrightarrow G$
maps $H$ to $\{1\}$. So $f_H(\gamma)\circ f_{G/H}(\gamma):G\longrightarrow G$ maps $G$ to $\{1\}$. We now
apply Lemma~\ref{4.1.1}.
\end{proof}

\begin{prop}\label{4.1.3}
Let $G$ be a group and let $\gamma$ be an endomorphism. Then the identities of $\gamma$ form an ideal of
$\Z[x]$.
\end{prop}

\begin{proof}
Suppose that $g_1(x)$ and $g_2(x)$ are identities of $\gamma$. For all $t\in G$, we then have
$(g_1(\gamma)(t))\cdot((g_2(\gamma))(t))=1\cdot 1=1$. So $g_1(x)+g_2(x)$ is again an identity of $\gamma$. Now
let $f(x)\in \Z[x]$ and let $g(x)$ be an identity of $\gamma$. Define $H:=G$. Then $f_{G/H}(x):=f(x)$ is an
identity of the natural endomorphism $\gamma_{G/H}$ induced on $G/H$ and $f_H(x):=g(x)$ is an identity of
$\gamma$ on $H$. So we need only apply Lemma~\ref{4.1.2} to conclude that also $f(x)\cdot g(x)$ is an identity
of $\gamma$.
\end{proof}

\subsection{Existence and construction of identities}\label{sec:4.2}

\begin{prop}\label{4.2.1}
Let $\alpha$ be a fixed-point-free automorphism of a finite non-trivial group $G$. Then $\alpha$ has an
ordered identity $f(x)$ that is monic and irreducible, that has no powerful divisors of positive degree, and
that satisfies $f(0)=-1$ and $f(1)\ne 0$.
\end{prop}

\begin{proof}
We note that $|\alpha|>1$ and $|G|>2$. Define the monic polynomial $f(x):=x^{|\alpha|}+ |G|\cdot
x^{|\alpha|-1}-1$. For all $t\in G$, we have $t^{-1}\cdot \alpha^{|\alpha|-1}(t^{|G|})\cdot
\alpha^{|\alpha|}(t)=t^{-1}\cdot 1\cdot t=1$, so that $f(x)$ is an ordered identity of $\alpha$. Perron's
criterion shows that $f(x)$ is irreducible. Since $f(x)$ is not powerful, all powerful divisors of $f(x)$ are
constant.
\end{proof}

\begin{prop}\label{4.2.2}
Let $G$ be a polycyclic group, let $\gamma$ be an endomorphism, and let $G=G_1\geqslant \cdots \geqslant
G_k\geqslant G_{k+1}=\{1_G\}$ be a subnormal, $\gamma$-invariant series such that the factors $G_i/G_{i+1}$
are elementary-abelian or free-abelian. Let $f_i(x)$ be an identity of the natural endomorphism
$\gamma_i:G_i/G_{i+1}\longrightarrow G_i/G_{i+1}$ induced on the factor $G_i/G_{i+1}$. Then
$f(x):=f_1(x)\cdots f_k(x)$ is an identity of $\gamma$.
\end{prop}

\begin{proof}
If $k\leqslant 1$, there is nothing to prove. Otherwise, we induct on $k$ to show that $f_2(x)\cdots f_k(x)$
is an identity of $\gamma$ on $G_2$. Lemma~\ref{4.1.2}, with $H:=G_2$, then finishes the proof.
\end{proof}

\begin{rem}\label{4.2.3}
$G$ has such a series since it is polycyclic. We also note that each endomorphism $\gamma_i$ has such a
(monic) identity by the theorem of Cayley--Hamilton.
\end{rem}

\begin{lem}\label{4.2.4}
A polynomial $f(x)\in \Z[x]$ is an identity of a fixed-point-free automorphism of a finite non-trivial group
if and only if $f(x)$ is not of the form $\pm t^m\cdot(t-1)^n$, where $m,n\in \Z_{\geqslant 0}$.
\end{lem}

\begin{proof}
Suppose first that $f(x)$ is not of the form $\pm x^m\cdot(x-1)^n$. Then $f(x)$ has an irreducible divisor
$s(x)$ that is not of that form either. If $s(x)$ is a prime $p$, then we may select any fixed-point-free
automorphism of any non-trivial elementary-abelian group of exponent $p$. Else, we write $s(x)=a_0+a_1\cdot
x+\cdots+a_d\cdot x^d$, with $a_0\cdot a_d\cdot s(1)\ne 0$ and $d\geqslant 1$. Let $q$ be any prime not
dividing $a_0\cdot a_d\cdot s(1)$. Then $s(x)$ is an identity of the fixed-point-free automorphism $a_d\cdot
C:\Z^d_q\longrightarrow \Z^d_q$, where $C$ is the companion operator of the monic polynomial $a^{-1}_d\cdot
s(x)$ in characteristic $q$. According to Proposition~\ref{4.1.3}, also $f(x)$ is an identity of the
automorphism. Now suppose that $\pm x^m\cdot(x-1)^n$ is an identity of a f.p.f. automorphism $\gamma$ of a
finite group $G$. Then $G$ is solvable according to Rowley's theorem~\ref{2.0.6}. Let $(G_i)_i$ be a series as
in Proposition~\ref{4.2.2}, which exists by Remark~\ref{4.2.3}. Then $\pm x^m\cdot(x-1)^n$ is an identity of
the endomorphism $\gamma_i$ that is induced on the factor $G_i/G_{i+1}$. According to Lemma~\ref{2.0.5}, each
$\gamma_i$ is a fixed-point-free automorphism. So the sections $G_i/G_{i+1}$ are trivial. So also $G$ is
trivial.
\end{proof}

The rather trivial construction in the proof of Lemma~\ref{4.2.4} can be refined as follows.

\begin{defi}\label{4.2.5}
Let $f(x)\in \Z[x]\setminus \Z$. We define $k(f(x))$ to be the maximal $k\in \Z_{\geqslant 1}$ such that there
exists a factorization $f(x)=f_1(x)\cdots f_k(x)$ of $f(x)$ with $f_1(x),\ldots,f_k(x)\in \Z[x]$ and roots
$\alpha,\beta \in \overline{\Q}$ with $f_1(\beta)=0$ and $f_1(\alpha)=f_2(\alpha \cdot \beta)=f_3(\alpha \cdot
\beta^2)=\cdots=f_k(\alpha \cdot \beta^{k-1})=0.$
\end{defi}

\begin{prop}\label{4.2.6}
Consider a monic polynomial $f(x)\in \Z[x]\setminus \Z$ with $f(0)=\pm 1$ and $f(1)\ne 0$. Let $p$ be any
prime not dividing $f(1)$. Then $f(x)$ is an identity of a fixed-point-free automorphism of a finite nilpotent
$p$-group $N$ of class exactly $k(f(x))$.
\end{prop}

Set $k:=k(f(x))$ and let $f(x)=f_1(x)\cdots f_k(x)$ be a decomposition as in Definition~\ref{4.2.5} with
corresponding roots $\alpha$ and $\beta$. According to Proposition~\ref{4.1.3}, may assume the $f_i(x)$ to be
square-free. We divide the proof into $3$ steps.

\begin{claim}\label{4.2.7}
There is a finite-dimensional, $k$-step nilpotent Lie algebra $L$ over $\Q$ with an automorphism
$A:L\longrightarrow L$ such that, for all $1\leqslant i\leqslant k$, the polynomial $f_i(x)$ is an identity of
the natural automorphism $A_i:\Gamma_i(L)/\Gamma_{i+1}(L)\longrightarrow \Gamma_i(L)/\Gamma_{i+1}(L)$ induced
on $\Gamma_i(L)/\Gamma_{i+1}(L)$.
\end{claim}

\begin{proof}
Let $F$ be the free $k$-step nilpotent Lie algebra over $\Q$ on the generators
$x_1,\ldots,x_d,\allowbreak y_1,\ldots,y_d$. For each $1\leqslant i\leqslant k$, we define $V_i$ to be the subspace
spanned by the left-normed Lie monomials of length exactly $i$ in the generators. Then $F=V_1\oplus \cdots
\oplus V_k$, as a vector space. Let $C$ be the companion operator of $f_1(x)$, which is semi-simple since
$f_1(x)$ is square-free. Let $A:V_1\longrightarrow V_1$ be the linear map with matrix $C\oplus C$ with respect
to the ordered basis $x_1,\ldots,y_d$. We may then extend this map (uniquely) to an automorphism
$A:F\longrightarrow F$ of $F$. For all $1\leqslant i\leqslant k$, we have $A(V_i)\subseteq V_i$, and therefore
$f(A)(V_i)\subseteq V_i$. We now define the ideal $I$ of $F$ that is generated by the set $f_1(A)(V_1)\cup
\ldots \cup f_k(A)(V_k)$ of homogeneous elements. Then $A(I)\subseteq I$ and $A^{-1}(I)\subseteq I$. Let
$\overline{A}:L\longrightarrow L$ be the natural automorphism induced on the quotient Lie algebra $L:=F/I$.
For each $1\leqslant i\leqslant k$, we clearly have $f_i(\overline{A})(\Gamma_i(L))\subseteq \Gamma_{i+1}(L)$.

We now suppose that $c(L)<k$ (and we will derive a contradiction). Then every left-normed Lie-monomial
$[u_1,\ldots,u_k]$ of length $k$ in the elements of $V_1$ is in the linear span of the sets
$[f_1(A)(V_1),V_1,\ldots,V_1],\ldots,[f_{k-1}(A)(V_{k-1}),V_1],f_k(A)(V_k).$ The same is true if we extend the
scalars from $\Q$ to $\C$. Let us now select an eigenbasis of $A$ on $V_1$, say
$v_1,\ldots,v_d,w_1,\ldots,w_d$. We may assume $A(v_1)=\alpha \cdot v_1$ and $A(w_1)=\beta \cdot w_1$. Then
the \emph{non-zero} left-normed monomial $[v_1,w_1,\ldots,w_1]$ of length $k$ is in the span of the
left-normed monomials
\[
[f_i(A)([u_1,\ldots,u_i]),u_{i+1},\ldots,u_k]=f_i(\lambda_1\cdots \lambda_i)\cdot[u_1,\ldots u_k],
\]
where $i$ runs over $\{1,,\ldots,k\}$ and where $u_1,\ldots,u_k$ are drawn from $\{v_1,\ldots,w_d\}$ with
corresponding eigenvalues $\lambda_1,\ldots,\lambda_k$. We may assume that, for each of these terms, the
multiplicity of $v_1$ and the multiplicity of $w_1$ in the sequence $u_1,\ldots,u_k$ is $1$ resp. $k-1$. Then
we have either $i\geqslant 2$ and $[u_1,\ldots,u_i]=[w_1,\ldots,w_1]=0$, or $f_i(\lambda_1\cdots
\lambda_i)=f_i(\alpha \cdot \beta^{i-1})=0$. So all terms vanish and $[v_1,w_1,\ldots,w_1]=0_F$ in the
\emph{free} Lie algebra $F$. This is our contradiction.
\end{proof}

\begin{claim}\label{4.2.8}
There is a finitely-generated, torsion-free, $k$-step nilpotent group $M$ with a fixed-point-free automorphism
$\widetilde{f}:M\longrightarrow M$ with identity $f(x)$.
\end{claim}

\begin{proof}
The Mal'cev correspondence (introduced in \cite{Malcev1949}) yields a group $(N,\cdot)$ and a fixed-point-free
automorphism $\widetilde{A}:N\longrightarrow N$ of $N$ such that the map $\widetilde{A}_i:
\Gamma_i(N)/\Gamma_{i+1}(N)\longrightarrow \Gamma_i(N)/\Gamma_{i+1}(N)$ induced on the $i$'th lower central
factor of $N$ coincides with the map $\overline{A}_i:\Gamma_i(L)/\Gamma_{i+1}(L)\longrightarrow
\Gamma_i(L)/\Gamma_{i+1}(L)$ on the $i$'th lower central factor of $L$. So $f_i(x)$ is an identity of
$\widetilde{A}_i$. Lemma~\ref{4.1.2} then implies that the polynomial $f(x)=f_1(x)\cdots f_k(x)$ is an
identity of $\widetilde{A}$. Since $f(1)=\pm 1$, we can select a full subgroup $M$ of $N$ that is invariant
under $\widetilde{A}$ and $\widetilde{A}^{-1}$ (using, for example, theorem $6.1$ of \cite{Dere}). The class
of $M$ is equal to the class of $N$ and $f(x)$ is an identity of the fixed-point-free automorphism
$\gamma:M\longrightarrow M$ that is obtained by restriction.
\end{proof}

\begin{claim}\label{4.2.9}
The polynomial $f(x)$ is an identity of a fixed-point-free automorphism of a finite nilpotent $p$-group $N$ of
class exactly $k$.
\end{claim}

\begin{proof}
According to Gruenberg's theorem~\ref{2.0.8}, $M$ is residually-(a finite $p$-group). So there exists a
characteristic subgroup $C$ of $M$ of $p$-power index such that $N:=M/C$ has class $k(f(x))$. Then $f(x)$ is
an identity of the naturally induced automorphism on $N$. Since $\gcd(|N|,f(1))=1$, the automorphism is
fixed-point-free.
\end{proof}

\section{Proof of Theorem~\ref{1.2.2}}\label{sec:5}

We continue to use the notation $f(x)=x^{m_0}\cdot(x-1)^{m_1}\cdot \overline{f}(x)$ of Definition~\ref{3.2.1}.

\begin{lem}\label{5.0.1}
Let $G$ be a finite, solvable group with a fixed-point-free automorphism $\alpha$ and let $f(x)$ be a non-zero
identity of the automorphism. Then $\overline{f}(x)$ is an identity of $\alpha$ on every $\alpha$-invariant,
elementary-abelian section of $G$.
\end{lem}

\begin{proof}
Consider an $\alpha$-invariant, elementary-abelian section $S$ of $G$. Then the induced automorphism
$\overline{\alpha}:S\longrightarrow S$ is fixed-point-free by Lemma~\ref{2.0.5}. So the maps $S\longrightarrow
S:s\mapsto \overline{\alpha}(s)$ and $S\longrightarrow S:s\mapsto s^{-1}\cdot \overline{\alpha}(s)$ are
surjective and we see that $(\overline{f}(\overline{\alpha}))(S)=(f(\overline{\alpha}))(S)= \{1\}$.
\end{proof}

\begin{thm}\label{5.0.2}
Let $G$ be a finite, solvable group with a fixed-point-free automorphism $\alpha$ and let $f(x)$ be an
identity of $\alpha$. Suppose that all powerful divisors of $\overline{f}(x)$ are constant. Then the first
lower Fitting group $\underline{F}_1(G)$ of $G$ is a $\pi(\iota_1(f(x)))$-group.
\end{thm}

\begin{proof}
We assume that $G$ is a counter-example of minimal order and we will derive a contradiction. If the
Fitting-height $h$ of $G$ is $1$ then $\underline{F}_1(G)=1$ is a $\pi(\iota_1)$-group (a contradiction). If
$h$ is at least three, then $\underline{F}_1(G)/\underline{F}_{h-1}(G)$ and
$\underline{F}_2(G)/\underline{F}_h(G)$ are $\pi(\iota_1)$-groups by the minimality of $G$. So, in this case,
also the extension $\underline{F}_1(G)$ is a $\pi(\iota_1)$-group (a contradiction). So we may assume that $G$
has Fitting-height exactly $2$. Let $R:= \underline{F}_1(G)$ and $F:=\overline{F}_1(G)$ be the lower resp.
upper Fitting group of $G$. Then we have the characteristic series $1\leqslant R\leqslant F\leqslant G$ of
$G$. If $R$ has two characteristic subgroups $A$ and $B$ with trivial intersection, then the minimality of $G$
forces $R/A$ and $R/B$ to be $\pi(\iota_1)$-groups, so that also $R/(A\cap B)=R/1$ is an $\pi(\iota_1)$-group
(a contradiction). So we may assume that $R$ is a $q$-group, for some prime $q$. Since $F$ is nilpotent, it is
the direct sum of a characteristic subgroup $M$ without $p$-torsion and a characteristic $q$-subgroup $N$,
naturally containing $R$. If $M\ne 1$, then the minimality of $G$ implies that $R/(R\cap M)=R/1$ is a
$\pi(\iota_1)$-group (a contradiction). So we may further assume that $F$ is a $q$-group.

We now consider the characteristic quotient $G/F$ and select an arbitrary, proper characteristic subgroup
$S/F$ of $G/F$. Then $S$ is a proper, characteristic subgroup of $G$. The minimality of $G$ implies that the
residual $R_0:=\Gamma_{\infty}(S)\subseteq R$ is a $\pi(\iota_1)$-group, characteristic in $G$. If $R_0\ne 1$,
then the minimality of $G$ implies that also $R/R_0$ is an $\pi(\iota_1)$-group, so that the extension $R$ is
an $\pi(\iota_1)$-group (a contradiction). So $R_0=1$ and $S$ is a normal, nilpotent subgroup of $G$ and
therefore contained in $F$. We conclude that $G/F$ is characteristically-simple. Since $G$ is solvable but not
nilpotent, we conclude that $G/F$ is an elementary-abelian $p$-group, for some prime $p$ unequal to $q$, say:
$G/F\cong(\Z^l_p,+)$. Lemma~\ref{2.0.5} yields a Hall-$p$-subgroup $P$ of $G$ that is $\alpha$-invariant and
such that $G=F\rtimes P.$

Now suppose that $R<F$. We then consider the characteristic quotient $G/R$. By construction, it is nilpotent,
so that $RP/R$ is characteristic in $G/R$. Its lift $RP=R\rtimes P$ is therefore characteristic in $G$. The
minimality of $G$ implies that the residual $R_1:= \Gamma_{\infty}(R\rtimes P)\subseteq R$ is a
$\pi(\iota_1)$-group, characteristic in $G$. If $R_1\ne 1$, then the minimality of $G$ implies that also
$R/R_1$ is a $\pi(\iota_1)$-group, so that the extension $R$ is a $\pi(\iota_1)$-group (a contradiction). So
$R_1=1$, $R\rtimes P$ is nilpotent, and $G$ is nilpotent (a contradiction). So $R=F$. The minimality of $G$
again implies that $R$ is elementary-abelian, say isomorphic to $(\Z^k_q,+)$ for some $k\in \Z_{\geqslant 1}$.

By embedding the group $(R,\cdot)$ into the additive group $(V,+)$ of the $k$-dimensional vector space $V$
over $\overline{\F}_q$, we naturally obtain groups of automorphisms $T,\langle A\rangle \subseteq \GL(V)$ with
the following properties (cf. Lemma~\ref{5.0.1}).
\begin{itemize}
\item[(a)]
$T\cong P$ and $T$ acts fixed-point-freely on $V$ by evaluation.
\item[(b)]
$A\in \mathrm{N}_{\GL(V)}(T)$ and $A$ acts fixed-point-freely on $T$ by conjugation within $\GL(V)$.
\item[(c)]
$\overline{f}(A)=0_V$ but $\iota_1\cdot V\ne 0_V$.
\end{itemize}

So we may consider the decomposition $V=\bigoplus_{\chi \in T^{\ast}}V_{\chi}$ of $V$ into its character
spaces $V_{\chi}=\{v\in V|\forall t\in T:t(v)=\chi(t)\cdot v\}$ w.r.t. the torus $T$. Since $\langle A\rangle$
normalizes $T$, each natural power $A^n$ of $A$ can be used to define a new character by the rule $(A^n\ast
\chi)(t):= \chi(A^{-n}\circ t\circ A^n).$ Then the group $\langle A\rangle$ acts on the set of character
spaces by the rule $A^n(V_{\chi})=V_{A^n\ast \chi}.$ Let us now fix an arbitrary, non-zero character space
$V_{\chi}$ corresponding to the character $\chi \in T^{\ast}$. Let $n$ be the cardinality of the orbit of
$V_{\chi}$ under the action of $\langle A\rangle$.

We claim that $n\geqslant 2$. Suppose first that, for every $t\in T$, the action of $t$ coincides with the
action of $A^{-1}\circ t\circ A$ on $V_{\chi}$. Since $A$ acts fixed-point-freely on $T$, we then see that
$\{[t,A]|t\in T\}=T$ acts trivially on all of $V_{\chi}\ne 0$, which contradicts our assumption $V^T=
\{0_V\}$. So there exists some $t\in T$ such that $A^{-1}\circ t\circ A$ and $A$ assume different values on
some $v\in V_{\chi}$. If $n=1$, then $f(V_{\chi})=V_{\chi}$, so that $(A^{-1}\circ t\circ
A)(v)=A^{-1}(t(A(v)))= A^{-1}(\chi(t)\cdot A(v))=\chi(t)\cdot v=t(v)$ (a contradiction).

Let us verify that $\iota_1$ vanishes on the linear subspace $W_{\chi}:=V_{\chi}\oplus A(V_{\chi})\oplus
\cdots \oplus A^{n-1}(V_{\chi})$ of $V$. Let $v$ be any vector of $V_{\chi}$. We first consider the
$n$-periodic decomposition $\overline{f}(x)=\overline{f}_{n,0}(x)+\cdots \overline{f}_{n,n-1}(x)$ of
$\overline{f}(x)$ of Definition~\ref{3.1.1}. By evaluating this decomposition in $A$, we obtain the equality
of $V$-endomorphisms: $\overline{f}(A)=\overline{f}_{n,0}(A)+\overline{f}_{n,1}(A)+
\cdots+\overline{f}_{n,n-1}(A).$ By further evaluating this decomposition in the distinguished vector $v\in
V_{\chi}$, we obtain the decomposition:
$(\overline{f}(A))(v)=(\overline{f}_{n,0}(A))(v)+\cdots+(\overline{f}_{n,n-1}(A))(v).$ By assumption, the map
$f(A)$ vanishes on all of $V$, so that the left-hand side of this decomposition is the trivial vector, $0_V$.
By construction, we also have the inclusions of vector spaces $\overline{f}_{n,0}(A)(V_{\chi})\subseteq
V_{\chi},\ldots,\overline{f}_{n,n-1}(A)(V_{\chi})\subseteq A^{n-1}(V_{\chi})$. Since the character-spaces
$V_{\chi}$, $A(V_{\chi}),\ldots,A^{n-1}(V_{\chi})$ are are distinct, they are linearly independent. So we may
conclude that all of the terms in this decomposition vanish:
\begin{equation}\label{eq:5.1}
(\overline{f}_{n,0}(A))(v)=\cdots=(\overline{f}_{n,n-1}(A))(v)=\{0_V\}.
\end{equation}
As mentioned in Remark~\ref{3.1.2}, there exist polynomials $g_0(x),\ldots,g_{n-1}(x)\in \Z[x]$ such that
$\iota_1=\iota_1(f(x))=g_0(x)\cdot \overline{f}_{n,0}(x)+\cdots+g_{n-1}(x)\cdot \overline{f}_{n,n-1}(x).$ By
evaluating in $A$, we obtain the decomposition of $V$-endomorphisms
\begin{equation}\label{eq:5.2}
\iota_1\cdot \mathbb{I}_V=g_0(A)\circ \overline{f}_{n,0}(A)+\cdots+g_{n-1}(A)\circ \overline{f}_{n,n-1}(A).
\end{equation}
By further evaluating this expression in the distinguished vector $v\in V_{\chi}$, we obtain
\[
\begin{aligned}
\iota_1\cdot v
 &\overset{\eqref{eq:5.2}}{=}
 (g_0(A)\circ\overline f_{n,0}(A))(v)+\cdots+
 (g_{n-1}(A)\circ\overline f_{n,n-1}(A))(v)\\
 &=g_0(A)(\overline f_{n,0}(A)(v))+\cdots+g_{n-1}(A)(\overline f_{n,n-1}(A)(v))\\
 &\overset{\eqref{eq:5.1}}{=}g_0(A)(0)+\cdots+g_{u-1}(A)(0)\\
 &=0_V.
\end{aligned}
\]
We conclude that the map $v\mapsto \iota_1\cdot v$ vanishes on each character space of $V$ and therefore on
all of $V$. But this contradicts our assumption (c).
\end{proof}

\begin{thm}\label{5.0.3}
Let $G$ be a finite group with a fixed-point-free automorphism $\alpha$ and let $f(x)$ be an identity of
$\alpha$. Suppose that all powerful divisors of $\overline{f}(x)$ are constant and suppose that
$\gcd(|G|,\iota_1(f(x)))=1.$ Then $G$ is nilpotent.
\end{thm}

\begin{proof}
We closely follow the strategy of \cite{Thompson1959a}. We suppose that, for a fixed $f(x)$, the conclusion is
false and we will derive a contradiction. Let $(G,\alpha,f(x))$ be a counter-example of minimal order $|G|$.
We note that every proper, $\alpha$-invariant section of $G$ is then nilpotent. If $|G|$ is the power of a
single prime, then $G$ is nilpotent (a contradiction). So some odd prime $p$ divides $|G|$. Lemma~\ref{2.0.5}
(b) gives us a Sylow $p$-subgroup satisfying $\alpha(P)=P$. We now distinguish between two cases and we show
that neither is possible. (a) Suppose first that $P$ has a non-trivial, normal, $\alpha$-invariant subgroup
$H$ with $\mathrm{N}_G(H)=G$. Then the section $G/H$ is proper and $\alpha$-invariant. So $G$ is the extension
of the nilpotent group $H$ by the nilpotent group $G/H$. The group $G$ is therefore solvable. Since
$\gcd(|G|,\iota_1(f(x)))=1$, Theorem~\ref{5.0.2} now implies that $G$ is nilpotent (a contradiction). (b)
Suppose next that every non-trivial, normal, $\alpha$-invariant subgroup $H$ of $P$ has a normalizer
$\mathrm{N}_G(H)$ that is properly contained in $G$. As a proper, $\alpha$-invariant section of $G$, this
$\mathrm{N}_G(H)$ is nilpotent. So Thompson's $p$-complement theorem~\ref{2.0.7} gives us a normal complement
$K$ to $P$ in $G$: $G\cong K\rtimes P$. As a proper, $\alpha$-invariant section of $G$, this $K$ is nilpotent.
As an extension of a nilpotent group $K$ by a nilpotent group $P$, our group $G$ is solvable.
Theorem~\ref{5.0.2} once more implies that $G$ is nilpotent (a contradiction).
\end{proof}

This establishes Theorem~\ref{1.2.2}, since every powerful divisor of $\overline{f}(x)$ is a divisor of
$f(x)$.

\section{Auxiliary theorems for Lie rings}\label{sec:6}

\subsection{An embedding theorem}\label{sec:6.1}

In this subsection alone, we will consider Lie rings with coefficients in a commutative ring $R$ that need not
be $\Z$. We refer to these as $R$-Lie rings.

\begin{defi}\label{6.1.1}
Let $n\in \Z$. We define the $n$-torsion ideal $T_n(L)$ of a ring $L$ by $T_n(L):=\{v\in L|\exists k\in
\Z_{\geqslant 1}:n^k\cdot v=0_L\}$. We say that $L$ \emph{has no} $n$\emph{-torsion} if $T_n(L)= \{0_L\}$.
\end{defi}

It is clear that the $n$-torsion ideal of a ring is an ideal of the ring that is invariant under all
endomorphisms. By definition, the additive group $(T_n(L),+)$ is a $\pi(n)$-group.

\begin{thm}\label{6.1.2}
Consider a $\Z$-Lie ring $L$, with an endomorophism $\gamma:L\longrightarrow L$, and a polynomial $f(x)\in
\Z[x]$ such that $f(\gamma)=0_L$. Suppose that $L$ has no $(\Discr_{\ast}(f(x))\cdot \Prod(f(x)))$-torsion.
Then there is an embedding of the $\Z$-Lie ring $\Discr_{\ast}(f(x))\cdot L$ into a $\Z$-Lie ring that is
graded by the semi-group $(\overline{\Q},\cdot)$ and supported by the root set $X$ of $f(x)$ in
$\overline{\Q}$.
\end{thm}

Let us use the abbreviations $\delta:=\Discr_{\ast}(f(x))$ and $\pi:=\Prod(f(x))$. We may assume that $f(x)$
is not a constant polynomial, since otherwise we trivially have $\delta \cdot L=0_L$, and there is nothing to
prove. Let $\lambda_1,\ldots,\lambda_l$ be the distinct roots with respective multiplicities $m_1,\ldots,m_l$
and let $a$ be the leading coefficient, so that $f(x)=a\cdot \prod_{1\leqslant i\leqslant
l}(x-\lambda_i)^{m_i}.$ Let $R:=\Z[\lambda_1,\ldots,\lambda_l]$ be the ring generated by the roots and let
$\F$ be the field of fractions of $R$. We next introduce a new Lie ring $\widetilde{L}:=R\otimes_{\Z}L$ with
coefficients in $R$. This $R$-Lie ring $\widetilde{L}$ naturally admits the $R$-Lie ring endomorphism
$\widetilde{\gamma}:\widetilde{L}\longrightarrow \widetilde{L}:\sum_ja_j\otimes v_j\longmapsto
\sum_ja_j\otimes(\gamma(v_j))$ and this map inherits the property $f(\widetilde{\gamma})=0_{\widetilde{L}}$.
Let $T$ be the $(\delta \cdot \pi)$-torsion ideal of $\widetilde{L}$. For each $\lambda \in R$, we define the
corresponding $R$-submodule of $\widetilde{L}$:
\[
E_{\lambda}:=\{v\in \widetilde{L}|\exists n\in \Z_{\geqslant 1}:(\widetilde{\gamma}-\lambda \cdot \mathbb{I}_{\widetilde{L}})^nv\in T\}.
\]

In what follows, we will repeatedly use Remark~\ref{3.1.6} without making this explicit.

\begin{claim}\label{6.1.3}
Let $\lambda,\mu \in R$, $m_{\lambda},m_{\mu}\in \Z_{\geqslant 1}$, and $v,w\in \widetilde{L}$. If
$(\widetilde{\gamma}-\lambda \cdot \mathbb{I}_{\widetilde{L}})^{m_{\lambda}}(v)=0_{\widetilde{L}}$ and
$(\widetilde{\gamma}-\mu \cdot \mathbb{I}_{\widetilde{L}})^{m_{\mu}}(w)=0_{\widetilde{L}}$, then
$(\widetilde{\gamma}-(\lambda \cdot \mu)\cdot
\mathbb{I}_{\widetilde{L}})^{(m_{\lambda}+m_{\mu})}([v,w])=0_{\widetilde{L}}$.
\end{claim}

\begin{proof}
We can prove, with Pascal's binomial identity and by induction on $m\in \Z_{\geqslant 1}$, that
$(\widetilde{\gamma}-(\lambda \cdot \mu)\cdot \mathbb{I}_{\widetilde{L}})^m([v,w])=\sum_{0\leqslant i\leqslant
m}\binom{m}{i}\cdot[\lambda^{m-i}\cdot(\widetilde{\gamma}-\lambda \cdot
\mathbb{I}_{\widetilde{L}})^i(v),\widetilde{\gamma}^i\circ(\widetilde{\gamma}-\mu \cdot
\mathbb{I}_{\widetilde{L}})^{m-i}(w)].$ So we need only specialize $m$ to $m_{\lambda}+m_{\mu}$.
\end{proof}

\begin{claim}\label{6.1.4}
Let $\lambda,\mu \in R$. Then we have (a) $T\subseteq \bigcap_{\nu \in R}E_{\nu}$, (b)
$[E_{\lambda},E_{\mu}]\subseteq E_{\lambda \cdot \mu},$ and (c) If $f(\lambda)=0=f(\mu)$ and $f(\lambda \cdot
\mu)\ne 0$, then $[E_{\lambda},E_{\mu}]\subseteq T$.
\end{claim}

\begin{proof}
(a) This property holds trivially. (b) Select arbitrary $v\in E_{\lambda}$ and $w\in E_{\mu}$. By definition,
there exist $n_1,n_2,n_3,n_4\in \Z_{\geqslant 1}$ such that $(\delta \cdot
\pi)^{n_1}\cdot(\widetilde{\gamma}-\lambda \cdot \mathbb{I}_{\widetilde{L}})^{n_2}v=0_{\widetilde{L}}=(\delta
\cdot \pi)^{n_3}(\widehat{\gamma}-\mu \cdot \mathbb{I}_{\widetilde{L}})^{n_4}w.$ By Claim~\ref{6.1.3}, we get:
$(\delta \cdot \pi)^{n_1+n_3}\cdot(\widetilde{\gamma}-\lambda \cdot \mu \cdot
\mathbb{I}_{\widetilde{L}})^{n_2+n_4}[v,w]=0_{\widetilde{L}}.$ (c) In view of the above, it suffices to show
that $E_{\lambda \cdot \mu}\subseteq T$. Let us abbreviate $\nu:=\lambda \cdot \mu$. Let $f(x)$ be given by
$f(x):=\sum_{0\leqslant i\leqslant d}a_i\cdot x^i\in \Z[x]$. Then $f(x)-f(\nu)=s(x)\cdot(x-\nu),$ where
$s(x)=\sum_{0\leqslant i\leqslant d}a_i\cdot(\sum_{j+k=i}x^j\cdot \nu^k)\in R[x]$. Select an arbitrary $v\in
E_{\nu}$. By definition, there exist $n_1,n_2\in \Z_{\geqslant 1}$ such that $(a\cdot \delta \cdot
\pi)^{n_1}\cdot(\widetilde{\gamma}-\nu \cdot \mathbb{I}_{\widetilde{L}})^{n_2}v=0_{\widetilde{L}}$. Then also
\begin{equation}\label{eq:6.1}
(\delta \cdot \pi)^{n_1}\cdot(f(\widetilde{\gamma})-f(\nu)\cdot \mathbb{I}_{\widetilde{L}})^{n_2}v=s(\widetilde{\gamma})^{n_2}\bigl((\delta \cdot \pi)^{n_1}\cdot(\widetilde{\gamma}-\nu \cdot \mathbb{I}_{\widetilde{L}})^{n_2}v\bigr)=0_{\widetilde{L}}.
\end{equation}
But, since $r(\widetilde{\gamma})$ vanishes on $\widetilde{L}$, we also have:
\begin{equation}\label{eq:6.2}
\begin{aligned}
(f(\widetilde\gamma)-f(\nu)\cdot\mathbb I_{\widetilde L})^{n_1}v
  &=\sum_{j+n_2=n_1}\binom{n_1}{j,n_2}\cdot f(\widetilde\gamma)^j
    \bigl((-f(\nu))^{n_2}\cdot v\bigr)\\
  &=(-f(\nu))^{n_1}\cdot v.
\end{aligned}
\end{equation}
By combining \eqref{eq:6.1} and \eqref{eq:6.2}, we obtain $(\delta \cdot \pi)^{n_2}\cdot r(\nu)^{n_1}\cdot
v=0_{\widetilde{L}}$. Since $f(\lambda)=0=f(\mu)$ and $f(\lambda \cdot \mu)\ne 0$, there exists a $n_3\in
\Z_{\geqslant 1}$ such that $(\delta \cdot \pi)^{n_2+n_3}\cdot v=0_{\widetilde{L}}$. So $v\in T$.
\end{proof}

\begin{claim}\label{6.1.5}
The $R$-module $K:=\sum_{\lambda \in X}E_{\lambda}$ is an $R$-Lie subring of $\widetilde{L}$.
\end{claim}

\begin{proof}
According to the Jacobi-identity and the bi-linearity of the Lie bracket, the Lie $R$-subalgebra of
$\widetilde{L}$ generated by the $R$-submodule $K$ of $\widetilde{L}$ is the $R$-span of left-normed words $w$
of the form $w:=[v_1,\ldots,v_n],$ where each $v_i$ is contained in some $K_{\mu_i}$ and each $\mu_i$ is
contained in $X$. So it suffices to show that for such a word $w$, we have $w\in K$. We may suppose that
$r(\mu_1)=\cdots=r(\mu_n)=0$, since otherwise $w$ is contained in the ideal $T$ and therefore in $K$. If
$\mu_1,\mu_1\cdot \mu_2,\ldots,\mu_1\cdots \mu_n:=\lambda$ are all roots of $f(x)$, then we need only apply
claim $1.b(n-1)$-times in order to conclude that $w\in K_{\lambda}$ and therefore $w\in K$. Else, there exists
an index $n_0\in \{1,\ldots,n-1\}$, such that $\mu_1\cdots \mu_{n_0}:=\mu$ is a root, but $\mu \cdot
\mu_{n_0+1}$ is not. Define $u:=[v_1,\ldots,v_{n_0}]$. By applying claim $1.b(n_0-1)$-times, we see that $u\in
K_{\mu}$. By applying claim 2.c, we see that $[u,v_{n_0+1}]\in T$. So also
$w=[[u,v_{n_0+1}],v_{n_0+2},\ldots,v_n]\in T\subseteq K$.
\end{proof}

\begin{claim}\label{6.1.6}
We have the inclusion $1\otimes(\delta \cdot L)\subseteq K$.
\end{claim}

\begin{proof}
For each $i\in \{1,\ldots,l\}$ and each $j\in \{0,\ldots,m_i-1\}$, we define the polynomial
$P_{i,j}(x):=f(x)/(x-\lambda_i)^{m_i-j}$ with coefficients in the ring $R:=\Z[\lambda_1,\ldots,\lambda_l]$.
Let us first show that, for each polynomial $P_{i,j}(x)$, there exists a coefficient $\theta_{i,j}\in R$ such
that
\begin{equation}\label{eq:6.3}
\Discr_{\ast}(f(x))=\sum_{1\leqslant i\leqslant l}\sum_{0\leqslant j\leqslant m_i-1}\theta_{i,j}\cdot P_{i,j}(x).
\end{equation}
The partial fraction decomposition of $a/f(x)$ is given by
\[
\frac{a}{f(x)}
 =\sum_{1\leqslant i\leqslant l}\sum_{0\leqslant j\leqslant m_i-1}
   \frac{1}{j!}\cdot\left(\left(\frac{a}{P_{i,0}(x)}\right)^{(j)}(\lambda_i)\right)
   \cdot(x-\lambda_i)^{j-m_i}.
\]
For each $i\in \{1,\ldots,l\}$, we define the auxiliary polynomial
\[
s_i(x):=\prod_{\substack{1\leqslant j\leqslant l\\j\ne i}}
(x-\lambda_j)\in R[x].
\]
Then the $j$'th derivative of $a/(P_{i,0}(x))$ is clearly of the form $b_{i,j}(x)/(s_i(x))^{2m}$, for some
explicitly computable $b_{i,j}(x)\in R[x]$. We see, in particular, that
\[
a=\sum_{1\leqslant i\leqslant l}\sum_{0\leqslant j\leqslant m_i-1}
 \frac{b_{i,j}(\lambda_i)}{j!\cdot(s_i(\lambda_i))^{2m}}\cdot P_{i,j}(x).
\]
After multiplying both sides of this equality by $\Discr_{\ast}(f(x))/a$, we see that
\[
\theta_{i,j}:=\frac{(m-1)!}{j!}\cdot
 \left((-1)^{m(l-1)}\cdot a^{2d^2}\cdot
 \prod_{\substack{1\leqslant k,n\leqslant l\\i\ne k\ne n\ne i}}
 (\lambda_k-\lambda_n)^m\right)\cdot b_{i,j}(\lambda_i)
\]
is a solution to \eqref{eq:6.3} in the ring of coefficients $R$, where $m:=\max \{m_1,\ldots,m_l\}$.

We now select an arbitrary $v\in L$. Corresponding with the $(i,j)$'th term of \ref{eq:6.3}, we define
$v_{i,j}:=P_{i,j}(\widetilde{\gamma})(1\otimes v)\in \widetilde{L}$. Then we observe that
$(\widetilde{\gamma}-\lambda_i\cdot \mathbb{I}_K)^{m_i-j}v_{i,j}=r(\widetilde{\gamma})(1\otimes v)=0_K$, so
that $v_{i,j}\in K_{\lambda_i}$. By evaluating the expression \eqref{eq:6.3} in $\widetilde{\gamma}$ and then
in $1\otimes v$, we see that
\[
1\otimes(\delta \cdot v)=\sum_{1\leqslant i\leqslant l}\sum_{0\leqslant j\leqslant m_i-1}\theta_{i,j}\cdot v_{i,j}\in \sum_{1\leqslant i\leqslant l}\sum_{0\leqslant j\leqslant m_i-1}R\cdot K_{\lambda_i}\subseteq K.
\]
So we may indeed conclude that $1\otimes(\delta \cdot L)\subseteq K$.
\end{proof}

We now note that $T$ is an ideal of $K$ contained in each $E_{\lambda}$. Let $\pi:K\longrightarrow K/T$ be the
natural projection. Then the $R$-Lie ring $\overline{K}:=\pi(K)$ is the sum $\sum_{\lambda \in
\overline{\Q}}\overline{K}_{\lambda}$ of the $R$-submodules
\[
\overline K_\lambda:=
\begin{cases}
 \pi(E_\lambda) & \text{if }\lambda\in X,\\
 \pi(T) & \text{if }\lambda\notin X.
\end{cases}
\]

\begin{claim}\label{6.1.7}
$\overline{K}=\bigoplus_{\lambda \in \overline{\Q}}\overline{K}_{\lambda}$ is a grading of the $R$-Lie ring
$\overline{K}$ by $(\overline{\Q}^{\times},\cdot)$ and its support is contained in $X$.
\end{claim}

\begin{proof}
In view of claims~\ref{6.1.4} and \ref{6.1.5}, we need only show that the decomposition
$\overline{K}=\sum_{\lambda \in \overline{\Q}}\overline{K}_{\lambda}$ of $\overline{K}$ into $R$-submodules is
direct. We select arbitrary $v_{\lambda_1}\in K_{\lambda_1},\ldots,v_{\lambda_l}\in K_{\lambda_l}$ such that
$v_{\lambda_1}+\cdots+v_{\lambda_l}\in T$ and we then need to show that $v_{\lambda_1},\ldots,v_{\lambda_l}\in
T$. So we select an arbitrary $i\in \{1,\ldots,l\}$ and will show that $v_{\lambda_i}\in T$. By definition,
there exist $n_1,n_2\in \Z_{\geqslant 1}$ such that for all $j\in \{1,\ldots,l\}$, we have:
\begin{equation}\label{eq:6.4}
(\delta \cdot \pi)^{n_2}\cdot(\widetilde{\gamma}-\lambda_j\cdot \mathbb{I}_K)^{n_1}v_{\lambda_j}=0_K
\end{equation}
Define the auxiliary polynomial $s(x)=\prod_{j\ne i}(x-\lambda_j)^{n_1}\in R[x]$. Then the theory of
resultants tells us that there exist polynomials $g(x),h(x)\in R[x]$ such that
\begin{equation}\label{eq:6.5}
g(x)\cdot s(x)+h(x)\cdot(x-\lambda_i)^{n_1}=\Res(s(x),(x-\lambda_i)^{n_1})=s(\lambda_i)^{n_1}.
\end{equation}
Since $T$ is invariant under $\widetilde{\gamma}$ and multiplication by elements of $R$, we see that $(\delta
\cdot \pi)^{n_2}\cdot s(\widetilde{\gamma})(v_{\lambda_i})=(\delta \cdot \pi)^{n_2}\cdot
s(\widetilde{\gamma})(v_{\lambda_1}+\cdots+v_{\lambda_l})\in T.$ So there exists a $n_3\in \Z_{\geqslant 1}$
such that
\begin{equation}\label{eq:6.6}
(\delta \cdot \pi)^{n_2+n_3}\cdot s(\widetilde{\gamma})(v_{\lambda_i})=0_K.
\end{equation}
By first evaluating \eqref{eq:6.5} in $\widetilde{\gamma}$ and then in $v_{\lambda_i}$, and by substituting
\eqref{eq:6.4} and \eqref{eq:6.6}, we obtain:
\[
\begin{aligned}
(\delta\cdot\pi)^{n_2+n_3}\cdot s(\lambda_i)^{n_1}v_{\lambda_i}
 &=(\delta\cdot\pi)^{n_2+n_3}\cdot g(\widetilde\gamma)(s(\widetilde\gamma)v_{\lambda_i})\\
 &\quad +(\delta\cdot\pi)^{n_2+n_3}\cdot
 h(\widetilde\gamma)((\widetilde\gamma-\lambda_i\cdot\mathbb I_K)^{n_1}v_{\lambda_i})\\
 &=0_K.
\end{aligned}
\]
Since the factor $s(\lambda_i)$ divides $\delta$ in the ring $R$, there exists some $n_4\in \Z_{\geqslant 1}$
such that also $(\delta \cdot \pi)^{n_2+n_3+n_2\cdot n_4}\cdot v_{\lambda_i}=0_K$. We may therefore conclude
that $v_{\lambda_i}\in T$.
\end{proof}

Any $R$-module $V$ is naturally also a $\Z$-module. Let us denote the corresponding $\Z$-module by $V_{\Z}$.
Then $K_{\Z}$ and $\overline{K}_{\Z}$ are naturally $\Z$-Lie rings. Moreover, the map $\pi_{\Z}:
K_{\Z}\longrightarrow \overline{K}_{\Z}:x\longmapsto \pi(x)$ is then a homomorphism of $\Z$-Lie rings.
According to Claim~\ref{6.1.6}, the map $\iota_{\Z}:\delta \cdot L\longrightarrow K_{\Z}:\delta \cdot
v\longmapsto 1\otimes(\delta \cdot v)$ is a well-defined homomorphism of $\Z$-Lie rings.

\begin{claim}\label{6.1.8}
The map $\pi_{\Z}\circ \iota_{\Z}:\delta \cdot L\longrightarrow \overline{K}_{\Z}:\delta \cdot v\longmapsto
\pi(1\otimes(\delta \cdot v))$ is an embedding of the $\Z$-Lie ring $\delta \cdot L$ into a $\Z$-Lie ring
$\overline{K}_{\Z}$ that is graded by $(\overline{\Q},\cdot)$ and supported by $X$.
\end{claim}

\begin{proof}
The map is a homomorphism of $\Z$-Lie rings since it is a composition of two such homomorphisms. The
homomorphism is injective, since $\operatorname{Ker}(\pi_{\Z}\circ \iota_{\Z})\subseteq T_{\delta \cdot
\pi}(L)$ and $T_{\delta \cdot \pi}(L)=\{0_L\}$ by assumption. Finally, the grading
$\overline{K}=\bigoplus_{\lambda \in \overline{\Q}}\overline{K}_{\lambda}$ of the $R$-Lie ring $\overline{K}$
in Claim~\ref{6.1.7} naturally gives rise to the grading $\overline{K}_{\Z}=\bigoplus_{\lambda \in
\overline{\Q}}(\overline{K}_{\lambda})_{\Z}$ of the $\Z$-Lie ring $\overline{K}_{\Z}$.
\end{proof}

This finishes the proof of Theorem~\ref{6.1.2}.

\subsection{A grading theorem}\label{sec:6.2}

\begin{thm}\label{6.2.1}
Let $X$ be a finite, arithmetically-free subset of the multiplicative group $(\F^{\times},\cdot)$ of a field
$\F$. Then every Lie ring $L$ that is graded by $(\F^{\times},\cdot)$ and supported by $X$ is nilpotent of
class at most $|X|^{2^{|X|}}.$
\end{thm}

We will prove the theorem by generalizing the proofs of well-known theorems of Kreknin--Kostrikin
\cite{KrekninKostrikin} and Khukhro \cite{Khukhro2008}, which correspond with a cyclic grading group of prime
order. We will also use Shalev's theorem~\ref{2.0.11}, which generalizes a theorem of Kreknin for cyclic
grading groups. Since Theorem~\ref{6.2.1} is a straightforward generalization of \cite{KrekninKostrikin}, we
present a compact proof.

We fix a finite, arithmetically-free subset $X$ of the multiplicative group $(\F^{\times},\cdot)$ of a field
$\F$. Clearly, $1\notin X$. We want to show that every Lie ring $L$ that is graded by $\F^{\times}$ and
supported by $X$ is nilpotent of class at most $|X|^{2^{|X|}}$. If $|X|=0$, then $L=0$, and we indeed have
$c(L)=0\leqslant 0^{2^0}$. If $|X|=1$, then $[L,L]\subseteq L_{G\setminus X}=\{0_L\}$, so that we have
$c(L)=1\leqslant 1^{2^1}$. We may therefore assume from now on that $x:=|X|\geqslant 2$. We begin by proving a
useful property of this arithmetically-free set $X$.

\begin{lem}\label{6.2.2}
For every $k\in \Z_{\geqslant 1}$, and elements $b_1,\ldots,b_k$ of $F\setminus \{1_F\}$, we can define the
set $S:=S(b_1,\ldots,b_k):=\bigcup_{0\leqslant i\leqslant k}\{b_{\sigma(1)}\cdots b_{\sigma(i)}|\sigma \in
\Sym(k)\}$ of partial products. Then: $|S|\geqslant k+1$, or $S$ contains a progression $g,g\cdot c,g\cdot
c^2,\ldots,g\cdot c^k$ with $c\in \{b_1,\ldots,b_k\}$.
\end{lem}

\begin{proof}
Let us prove the lemma by induction on $k\in \Z_{\geqslant 1}$. If $k=1$, then the statement is clearly true.
So we may suppose that $k>1$. We suppose that the second conclusion fails for $S$ and we will prove that
$|S|\geqslant k+1$. Define $S':=S(b_1,\ldots,b_{k-1})$ and $S:= S(b_1,\ldots,b_k)$. Then $S=S'\cup(S'\cdot
b_k)$ and the second conclusion also fails for $S'$. The induction hypothesis now tells us that $|S'|\geqslant
k$. If $S=S'\cdot b_k$, then $S=S'\cdot b_k=S'\cdot b^2_k= \cdots=S'\cdot b^k_k$. But then $S'$ would satisfy
condition $2$ (contradicting our assumption). So $S'\cdot b_k\ne S$, and therefore $|S|=|S'\cup(S'\cdot
b_k)|\geqslant 1+|S'|\geqslant 1+k$.
\end{proof}

\begin{lem}\label{6.2.3}
For all elements $a,b_1,\ldots,b_{|X|}$ of $X$, there exists a permutation $\sigma \in \Sym(|X|)$ and a
cut-off $k\in \{1,\ldots,|X|\}$ such that $a\cdot b_{\sigma(1)}\cdots b_{\sigma(k)}\in G\setminus X.$
\end{lem}

\begin{proof}
We suppose that the conclusion is false, and we will derive a contradiction. Consider
$S:=S(b_1,\ldots,b_{|X|})$. Then $a\cdot S\subseteq X$ and $|S|=|a\cdot S|\leqslant|X|$. Lemma~\ref{6.2.2} now
implies that $S$ contains a progression $b,b\cdot c,b\cdot c^2,\ldots,b\cdot c^{|X|}$ with $c\in X$. Then
$a\cdot S$, and therefore $X$, contains the progression $a\cdot b,(a\cdot b)\cdot c,(a\cdot b)\cdot
c^2,\ldots,(a\cdot b)\cdot c^{|X|}$ with $c\in X$. We conclude that $X$ is not arithmetically-free. This is a
contradiction.
\end{proof}

Let us use, for $u\in L$ and $k\in \Z_{\geqslant 0}$, the abbreviation $[u,{}_kL]:=[u,L,\ldots,L]$, where $L$
appears $k$ times on the right hand side. And let $(\Gamma_i(L))_{i\geqslant 1}$ and
$(\Delta_i(L))_{i\geqslant 0}$ be the lower and descending central series of $L$.

\begin{lem}\label{6.2.4}
Consider a Lie ring $L$. For every permutation $\sigma \in \Sym(X)$ and for all elements $v,w_1,\ldots,w_x\in
L$, we have $[v,w_1,\ldots,w_x]\in[v,w_{\sigma(1)},\ldots,w_{\sigma(x)}]+ \operatorname{id}_L([v,[L,L]]).$
\end{lem}

\begin{proof}
A repeated application of the Jacobi-identity shows that the vector
$[v,w_1,\ldots,w_x]-[v,w_{\sigma(1)},\ldots,w_{\sigma(x)}]$ is in the $\Z$-linear span of $[v,[L,L]]$,
$[v,[L,L],L]$, $\ldots,[v,[L,L],{}_xL],$ and therefore in the ideal of $L$ that is generated by $[v,[L,L]]$.
\end{proof}

\begin{prop}\label{6.2.5}
Consider a Lie ring $L$ that is graded by $(\F^{\times},\cdot)$ and supported by $X$.
\begin{itemize}
\item[(a)]
For every homogeneous ideal $I$ of $L$, we have $[I,{}_xL]\subseteq[I,[L,L]].$
\item[(b)]
For every homogeneous ideal $I$ of $L$ and $l\in \Z_{\geqslant 1}$, we have
$[I,{}_{lx}L]\subseteq[I,{}_l[L,L]].$
\item[(c)]
For every $l\in \Z_{\geqslant 1}$, we have the inclusion of ideals: $\Gamma_{2+lx}(L)\subseteq
\Gamma_{l+1}(\Delta_1(L)).$
\item[(d)]
For every $l\in \Z_{\geqslant 1}$, we have the inclusion of ideals $\Gamma_{1+(x^s-1)/(x-1)}(L)\subseteq
\Delta_s(L).$
\end{itemize}
\end{prop}

\begin{proof}
Let the gradings be given by $L=\bigoplus_{\lambda \in \F^{\times}}L_{\lambda}$ and $I=\bigoplus_{\lambda \in
\F^{\times}}I_{\lambda}$. (a) It suffices to show that every left-normed bracket of the form
$u:=[v_a,w_{b_1},\ldots,w_{b_x}],$ with $v_a\in I_a$ and $w_{b_1}\in L_{b_1},\ldots,v_{b_x}\in L_{b_x}$, is in
the ideal $[I,[L,L]]$. If any of the $a,b_1,\ldots,b_x$ are in $\F^{\times}\setminus X$, then we indeed have
$u=0\in[I,[L,L]]$. So we may assume that $a,b_1,\ldots,b_x\in X$. According to Lemma~\ref{6.2.3}, there exists
a permutation $\sigma \in \Sym(X)$ and a $k\in \{1,\ldots,x\}$, such that $c:=a\cdot b_{\sigma(1)}\cdots
b_{\sigma(k)}\in \F^{\times}\setminus X$, and therefore $L_c=\{0_L\}$. We then apply Lemma~\ref{6.2.4} to the
bracket $u$ and $\sigma$ in order to obtain
$u\in[[v_a,w_{\sigma(1)},\ldots,w_{\sigma(k)}],w_{\sigma(k+1)},\ldots,w_{\sigma(x)}]+[I,[L,L]].$ The grading
property implies that the bracket $[v_a,w_{\sigma(1)},\ldots,w_{\sigma(k)}]$ is contained in the homogeneous
component $L_c=\{0_L\}$ of $L$. So $u\in[0,w_{\sigma(k+1),\ldots,w_{\sigma(x)}}]+[I,[L,L]]=[I,[L,L]]$. (b) We
proceed by induction on $l\in \Z_{\geqslant 1}$. The base case $l=1$ corresponds with the above claim. So we
may assume that $l>1$. Then $J:=[I,{}_{(l-1)x}L]$ is a homogeneous ideal of $L$ and the first claim gives us
the inclusion $[J,{}_xL]\subseteq[J,[L,L]].$ The induction hypothesis for $l-1$ also gives the inclusion
$J\subseteq[I,{}_{l-1}[L,L]]$. These two combine to give
$[I,{}_{lx}L]\subseteq[[I,{}_{l-1}[L,L]],[L,L]]=[I,{}_l[L,L]].$ (c) We need only specialize the second claim
to the ideal $I:=[L,L]$ of $L$. (d) We proceed by induction on $s$. If $s=1$, then
$\Gamma_2(L)=[L,L]=\Delta_1(L)$. Now suppose that $s>1$. We then define $l:=(x^{s-1}-1)/(x-1)$ and we observe
that $(x^s-1)/(x-1)=1+xl$. According to the third claim, we have $\Gamma_{2+xl}(L)\subseteq
\Gamma_{1+l}(\Delta_1(L))$. We note that the Lie ring $\Delta_1(L)$ is naturally graded by
$(\F^{\times},\cdot)$ and supported by $X$. So the induction hypothesis for $s-1$ tells us that
$\Gamma_{l+1}(\Delta_1(L))\subseteq \Delta_{s-1}(\Delta_1(L))=\Delta_s(L)$. We therefore obtain
$\Gamma_{2+xl}(L)\subseteq \Gamma_{1+l}(\Delta_1(L))\subseteq \Delta_s(L).$
\end{proof}

\begin{proof}[Proof of Theorem~\ref{6.2.1}]
As above, we may suppose that $x:=|X|\geqslant 2$. Since $X$ is an arithmetically-free subset of
$\F^{\times}$, it does not contain $1$. According to Theorem~\ref{2.0.11}, the Lie ring $L$ is solvable of
derived length $\dl(L)\leqslant 2^x=:s$. Proposition~\ref{6.2.5} gives $\Gamma_{x^{2^x}+1}(L)\subseteq
\Gamma_{1+(x^s-1)/(x-1)}(K)\subseteq \Delta_s(L)\subseteq \Delta_{\dl(L)}(L)=\{0_L\}.$
\end{proof}

\Needspace{5\baselineskip}
\subsection{Lie rings with an endomorphism satisfying an identity}\label{sec:6.3}

\begin{thm}\label{6.3.1}
Let $L$ be a Lie ring, let $\gamma:L\longrightarrow L$ be an endomorphism, and let $f(x)$ be an identity of
$\gamma$. Suppose that the roots of $f(x)$ form an arithmetically-free subset $X$ of
$(\overline{\Q}^{\times},\cdot)$. Then the group $(\Gamma_{H(X,\overline{\Q}^{\times})+1}(L),+)$ and its
subgroup $(\Gamma_{|X|^{2^{|X|}}+1}(L),+)$ are $\pi(\iota_2)$-groups.
\end{thm}

\begin{proof}
We consider the $\iota_2$-torsion ideal $T_{\iota_2}$ of $L$. Then the Lie ring $M:=L/T_{\iota_2}$ has no
$\iota_2$-torsion and the induced endomorphism $\gamma:M\longrightarrow M$ satisfies $f(\gamma)=0_M$.
Theorem~\ref{6.1.2} embeds $\Discr_{\ast}(f(x))\cdot M$ into a Lie ring that is graded by
$(\overline{\Q}^{\times},\cdot)$ and supported by $X$. Theorem~\ref{2.0.10} shows that
$\Discr_{\ast}(f(x))\cdot M$ is nilpotent of class at most $H(X,\overline{\Q}^{\times})$, so that
$\Discr_{\ast}(f(x))^{H(X,\overline{\Q}^{\times})+1}\cdot
\Gamma_{H(X,\overline{\Q}^{\times})+1}(M)=\Gamma_{H(X,\overline{\Q}^{\times})+1}(\Discr_{\ast}(f(x))\cdot
M)=\{0\}$. Since $M$ has no $\iota_2$-torsion, we have $\Gamma_{H(X,\overline{\Q}^{\times})+1}(M)=\{0\}$, and
therefore $\Gamma_{H(X,\overline{\Q}^{\times})+1}(L)\subseteq T_{\iota_2}$. Theorem~\ref{6.2.1} states that
$H(X,\overline{\Q}^{\times})\leqslant|X|^{2^{|X|}}$, so that also $(\Gamma_{|X|^{2^{|X|}}+1}(L),+)$ is a
$\pi(\iota_2)$-group.
\end{proof}

\section{Proof of Theorem~\ref{1.2.4}}\label{sec:7}

\begin{thm}\label{7.0.1}
Let $G$ be a nilpotent group, let $\gamma$ be an endomorphism, and let $f(x)$ be an identity of the
endomorphism. Suppose that the roots of $f(x)$ form an arithmetically-free subset of
$(\overline{\Q}^{\times},\cdot)$. Then $\Gamma_{H(X,\overline{\Q}^{\times})+1}(G)$ and its subgroup
$\Gamma_{|X|^{2^{|X|}}+1}(G)$ are $\pi(\iota_2(f(x)))$-groups.
\end{thm}

\begin{proof}
We use standard properties of associated Lie rings. Cf. paragraph $3.2$ of \cite{Khukhro1993} for details. We
consider the Lie ring $L:=\bigoplus_{n\in \Z_{\geqslant 1}}\Gamma_n(G)/\Gamma_{n+1}(G)$ of $G$ that is
associated with the lower central series of $N$ and the natural Lie ring endomorphism $\gamma:L\longrightarrow
L$ that is induced on it. Then $\gamma$ satisfies $f(\gamma)=0_L$. Theorem~\ref{6.3.1} states that
$(\Gamma_{H(X,\overline{\Q}^{\times})+1}(L),+)$ is a $\pi(\iota_2(f(x)))$-group. Now note that, for all
$i\geqslant H(X,\overline{\Q}^{\times})+1$, the group $(\Gamma_i(G)/\Gamma_{i+1}(G),\cdot)$ is isomorphic to
the $\pi(\iota_2(f(x)))$-group $(\Gamma_i(L)/\Gamma_{i+1}(L),+)$. Since $G$ is nilpotent,
$\Gamma_{H(X,\overline{\Q}^{\times})+1}(G)$ is an extension of \emph{finitely}-many such groups and therefore
itself a $\pi(\iota_2(f(x)))$-group.
\end{proof}

Theorem~\ref{1.2.4} is a direct consequence of the above theorem.

\section{Proof of the corollaries}\label{sec:8}

\subsection{Proof of Corollary~\ref{1.0.1}}\label{sec:8.1}

\begin{defi}\label{8.1.1}
Let $f(x)\in \Z[x]$ be irreducible. If $f(x)\in \{\pm x,\pm(x-1)\}$, then we define $h(x):=1$ and
$\iota(f(x)):=1$. If $f(x)$ is a constant, then we define $h(x):=f(x)$ and $\iota(f(x)):=f(x)$. Else, we
consider the (necessarily unique) largest $m\in \Z_{\geqslant 1}$ such that $f(x)\in \Z[x^m]$. Then there is a
(necessarily unqiue) $h(x)\in \Z[x]$ such that $f(x)=h(x^m)$. We then define $\iota(f(x)):=h(1)\cdot
\iota_1(h(x))\cdot \iota_2(f(x))$.
\end{defi}

We continue to use the notation of the definition.

\begin{lem}\label{8.1.2}
Let $f(x)\in \Z[x]$ be irreducible. Then $\iota(f(x))\in \Z \setminus \{0\}$ and $h(x)$ has no powerful
divisors of positive degree.
\end{lem}

\begin{proof}
We may assume that $f(x)\notin \{\pm x,\pm(x-1)\}$ and that $f(x)$ is not constant, since otherwise the proof
is trivial. Then the polynomial $h(x)$ is again irreducible and, by construction, it is not powerful. So the
polynomial $h(x)$ has no powerful divisors of positive degree. Since $h(0)=f(0)\ne 0$,
Propositions~\ref{3.1.7} and \ref{3.2.2} then show that $\iota(f(x))\in \Z \setminus \{0\}$.
\end{proof}

\begin{lem}\label{8.1.3}
Let $G$ be a finite group with an automorphism $\alpha$ and let $f(x)$ be an identity of $\alpha$. Suppose
that $\gcd(|G|,f(1))=1$. Then $\alpha$ is fixed-point-free.
\end{lem}

\begin{proof}
By definition, there exist $f_1(x),\ldots,f_k(x)\in \Z[x]$ such that $f(x)=f_1(x)+ \cdots+f_k(x)$ and
$((f_1(\alpha))h)\cdots((f_k(\alpha))h)=1$ for all $h\in G$. Now suppose $\alpha(h)=h$. Then
$h^{f(1)}=h^{f_1(1)}\cdots h^{f_k(1)}=((f_1(\alpha))h)\cdots((f_k(\alpha))h)=1$. Since $\gcd(|G|,f(1))=1$, we
have $h=1$.
\end{proof}

\begin{cor}\label{8.1.4}
Let $f(x)\in \Z[x]$ be irreducible with root set $X$. Let $G$ be a finite group with a fixed-point-free
automorphism $\alpha$. Suppose that $f(x)$ is an ordered identity of $\alpha$ and that
$\gcd(|G|,\iota(f(x)))=1$. Then $G$ is nilpotent and $c(G)\leqslant H(X\setminus
\{0,1\},\overline{\Q}^{\times})\leqslant|X|^{2^{|X|}}.$
\end{cor}

\begin{proof}
If $f(x)\in \{\pm x\}$, then $G$ coincides with the kernel of $\alpha$. Since $\alpha$ is an automorphism, we
see that then $G=\{1\}$, so that there is nothing to prove. If $f(x)\in \{\pm(x-1)\}$, then $G$ coincides with
the fixed-points of $\alpha$. So we may again conclude that $G=\{1\}$. If $f(x)$ is constant, then the order
of every element of $G$ divides $f(x)$. Since $\gcd(|G|,f(x))= 1$, we may also in this case conclude that
$G=\{1\}$. So we may assume that $f(x)$ is none of the above. Then $X=X\setminus \{0,1\}$ is an
arithmetically-free subset of $(\overline{\Q}^{\times},\cdot)$, by Example~\ref{2.0.9}. Now let $f(x)=h(x^m)$
be the factorization of Definition~\ref{8.1.1}. Then $h(x)$ is an ordered identity of the automorphism
$\alpha^m$ and $h(x)$ has no powerful divisors of positive degree (Lemma~\ref{8.1.2}). Since
$\gcd(|G|,h(1))=1$, this automorphism $\alpha^m$ is fixed-point-free (Lemma~\ref{8.1.3}). So we may apply
Theorem~\ref{1.2.2} to conclude that $G$ is nilpotent. Theorem~\ref{1.2.4} gives the bound
$c(G)\leqslant|X|^{2^{|X|}}$ of Corollary~\ref{1.0.1}, while \ref{7.0.1} gives the more precise bounds
$c(G)\leqslant H(X,\overline{\Q}^{\times})\leqslant|X|^{2^{|X|}}$.
\end{proof}

\begin{rem}\label{8.1.5}
\emph{Almost-every} monic polynomial $f(x)\in \Z[x]$ in the standard models of
\cite{BKK,BK,BV,Eberhard,Rivin,VanDerWaerden} is irreducible and has a product-free root set $X$, so that
$H(X\setminus \{0,1\},\overline{\Q}^{\times})\leqslant 1$, and so that $G$ is \emph{abelian}. This bound
$c(G)\leqslant 1$ is optimal in view of Lemma~\ref{4.2.4}.
\end{rem}

This establishes Corollary~\ref{1.0.1}. By (implicitly) using the classification of the finite simple groups,
we obtain:

\begin{cor}\label{8.1.6}
Let $f(x)\in \Z[x]$ be irreducible with root set $X$. Let $G$ be a finite group with a fixed-point-free
automorphism $\alpha$ and suppose that $f(x)$ is an ordered identity of $\alpha$. Then $G$ is solvable and the
product $G=A\cdot B$ of a $\pi(\iota(f(x)))$-group $A$ and a nilpotent $\pi(\iota(f(x)))'$-group $B$ of class
at most $|X|^{2^{|X|}}$.
\end{cor}

\begin{proof}
The group is solvable according to Theorem~\ref{2.0.6}. So Lemma~\ref{2.0.5} gives us a Hall
$\pi(\iota(f(x)))$-subgroup $A$ and a Hall $\pi(\iota(f(x)))'$-subgroup $B$ of $G$ such that $\alpha(B)=B$.
Then $G=A\cdot B$ and we need only apply Corollary~\ref{8.1.4} to $B$, the restriction
$\alpha:B\longrightarrow B$, and $f(x)$.
\end{proof}

\begin{rem}\label{8.1.7}
This paper is the first in a series of three studying the structure of finite groups $G$ with a
fixed-point-free automorphism satisfying a given (non-zero) polynomial identity $f(x)$. In this paper, we have
assumed $f(x)$ to be \emph{irreducible} and we have obtained, among other things, the nilpotency of the group
modulo some prescribed torsion. In the third paper, joint with E.I. Khukhro, we consider the most general case
of non-zero identities $f(x)$. In view of Gross' examples \cite{Gross}, the conclusions there are necessarily
weaker than those of Corollary~\ref{1.0.1}: we obtain a rather large upper bound on the \emph{Fitting height}
of the group, in the spirit of Thompson's theorem \cite{Thompson1964} and Berger's theorem \cite{Berger}. The
second paper obtains a better bound on the Fitting height of $G$ (namely, the number of irreducible factors of
$f(x)$) if $f(x)$ belongs to the large family of \emph{Higman-solvable} polynomials. We refer to the
respective papers for more context and details.
\end{rem}

\subsection{Proof of Corollary~\ref{1.2.7}}\label{sec:8.2}

Let $\Kh $ be the invariant of Khukhro's Theorem~\ref{2.0.13} and let $n_0$ be the \emph{square-free radical}
of a given natural number $n$. We will repeatedly use the elementary properties of cyclotomic polynomials that
are listed in Section~\ref{sec:3}.

\begin{defi}\label{8.2.1}
Let $k\in \Z_{\geqslant 1}$ and $n\in \Z_{\geqslant 2}$. We define
\[
\mathrm B(k,n):=
\begin{cases}
 \max\{\varphi(n_0)^{2^{\varphi(n_0)}},\Kh(k,n_0)\}
    & \text{if $n_0$ is a prime},\\[2pt]
 \varphi(n_0)^{2^{\varphi(n_0)}} & \text{if $n_0$ is composite}.
\end{cases}
\]
\end{defi}

\begin{cor}\label{8.2.2}
Let $n\in \Z_{\geqslant 2}$ and $k\in \Z_{\geqslant 1}$. Let $G$ be a finite group on $k$ generators with an
$n$-cyclotomic automorphism $\alpha$. Then $G$ is nilpotent and $c(G)\leqslant \mathrm{B}(k,n)$.
\end{cor}

\begin{proof}
The automorphism $\beta:=\alpha^{n/n_0}$ is $n_0$-cyclotomic according to \eqref{eq:3.1}. Suppose first that
$n_0$ is a prime. Then $G$ is nilpotent according to Theorem~\ref{2.0.12} of Hughes--Thompson and Kegel. The
Hall $\pi(n_0)$-subgroup of $G$ has class at most $\mathrm{B}(k,n)$ according to Khukhro's
Theorem~\ref{2.0.13}. The Hall $\pi(n_0)'$-subgroup of $G$ has class at most $\mathrm{B}(k,n)$ according to
the theorem of Kreknin--Kostrikin \cite{KrekninKostrikin}. Instead of the latter, we may also combine
Theorem~\ref{1.2.4} with Proposition~\ref{3.3.5}. Either way, we conclude that $c(G)\leqslant
\mathrm{B}(k,n)$.

So we may assume that $n_0$ is composite. Since $\Phi_{n_0}(1)=1$, the automorphism $\beta$ is
fixed-point-free (Lemma~\ref{8.1.3}). Note that $\iota_1(\Phi_{n_0}(x))=1$ (Proposition~\ref{3.3.4}). So the
polynomial $\Phi_{n_0}(x)$ has no powerful divisors of positive degree (Proposition~\ref{3.1.4}). According to
Theorem~\ref{1.2.2}, $G$ is therefore nilpotent. It now suffices to show that the class of any Sylow
$p$-subgroup $P$ of $G$ is at most $l:=\varphi(n_0)^{2^{\varphi(n_0)}}$. Let $L$ be the Lie ring of such a
subgroup $P$ corresponding with the lower central series of $P$. Since the class of $L$ coincides with the
class of $P$, we need only show that $c(L)\leqslant l$. Let $\overline{\beta}$ be the automorphism that is
induced by $\beta$ on $L$. Then $\Phi_{n_0}(x)$ is an identity of $\overline{\beta}$. Since $n_0$ is
composite, it has a prime divisor $q$ that is distinct from $p$. We recall that
$\Phi_{n_0}(x)|\Phi_q(x^{n_0/q})$ (cf. \eqref{eq:3.2}). According to Proposition~\ref{4.1.3}, the polynomial
$\Phi_q(x^{n_0/q})$ is then an identity of $\overline{\beta}$. Since the additive group of $L$ is abelian, the
polynomial $\Phi_q(x)$ is then an identity of the automorphism $\overline{\beta}^{n_0/q}$. Moreover,
$\Phi_q(x)$ has an arithmetically-free set of roots (cf. Example~\ref{2.0.9}). So we may apply
Theorem~\ref{6.3.1} to the Lie ring $L$ with the automorphism $\overline{\beta}^{n_0/q}:L\longrightarrow L$
and identity $\Phi_q(x)$. We conclude that $(\Gamma_{\varphi(q)^{2^{\varphi(q)}}+1}(L),+)$ is a
$\pi(\iota_2(\Phi_q(x)))$-group. According to Proposition~\ref{3.3.5},
$(\Gamma_{\varphi(q)^{2^{\varphi(q)}}+1}(L),+)$ is therefore a $q$-group. But, by construction, we have
$|P|=|L|$. So $(L,+)$, and therefore $(\Gamma_{\varphi(q)^{2^{\varphi(q)}}+1}(L),+)$, is a $p$-group. Since
$p\ne q$, we conclude that $c(L)\leqslant \varphi(q)^{2^{\varphi(q)}}\leqslant l$. This finishes the proof.
\end{proof}

This establishes Corollary~\ref{1.2.7}. We can extend the result to compact groups as follows.

\begin{cor}\label{8.2.3}
Let $n\in \Z_{\geqslant 2}$ and $k\in \Z_{\geqslant 1}$. In a compact Hausdorff group $G$ with an
$n$-cyclotomic automorphism $\alpha$, every $k$-generated subgroup $H$ is nilpotent of class at most
$\mathrm{B}((2\cdot \varphi(n_0)-1)\cdot k,n)$.
\end{cor}

\begin{proof}
Let $H$ be abstractly generated by the $k$ elements $x_1,\ldots,x_k$. Then the $k_n:= (2\varphi(n)-1)k$
elements $\{\alpha^j(x_i)|1\leqslant i\leqslant k;-\varphi(n)<j<\varphi(n)\}$ generate an $\alpha$-invariant
subgroup $\widetilde{H}$ of $G$. According to Theorem~\ref{2.0.14} of Peter--Weyl, the group $G$, and
therefore its subgroup $\widetilde{H}$, is residually-(a linear group in characteristic $0$). Mal'cev's
Theorem~\ref{2.0.15} gives us a family $\mathcal{C}$ of normal subgroups of $\widetilde{H}$ of finite index
and with trivial intersection. Since $\widetilde{H}$ is finitely-generated, we may even assume these subgroups
to be characteristic in $\widetilde{H}$. For each $C\in \mathcal{C}$, we then have the naturally induced
$n$-cyclotomic automorphism $\alpha_C:\widetilde{H}/C\longrightarrow \widetilde{H}/C$ on the finite,
$k_n$-generated group $\widetilde{H}/C$. Corollary~\ref{8.2.2} now states that $\widetilde{H}/C$ is nilpotent
of class $c(\widetilde{H}/C)\leqslant \mathrm{B}(k_n,n)$. So $\Gamma_{\mathrm{B}(k_n,n)+1}(H)\subseteq
\Gamma_{\mathrm{B}(k_n,n)+1}(\widetilde{H})\subseteq \cap_{C\in \mathcal{C}}C=\{1_{\widetilde{H}}\}$. The
Hausdorff property implies that also the closure of $H$ in $G$ is nilpotent of class at most
$\mathrm{B}(k_n,n)$.
\end{proof}

\section*{Acknowledgements}

Dietrich Burde is supported by the FWF, the Austrian Science Foundation, with the 
grant-DOI 10.55776/PAT3410826.
Wolfgang A. Moens was supported by the FWF through grants J~3371 and P~30842.

\end{document}